\documentclass[10pt]{article}
\usepackage{amsmath, amssymb,amsthm,color,enumerate,comment,hyperref,graphicx,float,appendix,amssymb,enumitem,subcaption,booktabs,arydshln,mathrsfs}
\newcommand{\no}[1]{#1}
\renewcommand{\d}{\,\mathrm{d}}
\renewcommand{\no}[1]{}
\no{\usepackage{times}\usepackage[subscriptcorrection, slantedGreek, nofontinfo]{mtpro}
\renewcommand{\Delta}{\upDelta}}

\title{Stable Determination and Reconstruction of a Quasilinear Term in an Elliptic Equation\thanks{The work of M. Deng is supported by the Hong Kong PhD fellowship scheme, and that of B. Jin is
supported by Hong Kong RGC General Research Fund (Projects 14306423 and 14306824), ANR / RGC Joint Research
Scheme (A-CUHK402/24) and a start-up fund from The Chinese University of Hong Kong. The work of Y. Kian is supported by the French National Research Agency ANR and Hong Kong RGC Joint Research Scheme for the project IdiAnoDiff (grant ANR-24-CE40-7039).}}
\date{}
\author{Jason Choy\thanks{Department of Mathematics, The Chinese University of Hong Kong, Shatin, N.T., Hong Kong (\texttt{zhchoy@math.cuhk.edu.hk, mldeng@link.cuhk.edu.hk, b.jin@cuhk.edu.hk})}\and Maolin Deng\footnotemark[2] \and Bangti Jin\footnotemark[2] \and Yavar Kian\thanks{Universit\'e de Rouen Normandie, CNRS UMR 6085, Math\'ematiques, 76801 Saint-Etienne du Rouvray, France (\texttt{yavar.kian@univ-rouen.fr})}}

\date{\today}
\newtheorem{theorem}{Theorem}[section]
\newtheorem{proposition}{Proposition}[section]

\newtheorem{remark}{Remark}[section]

\newcommand{\R}{{\mathbb R}}

\numberwithin{equation}{section}

\renewcommand{\leq}{\leqslant}
\renewcommand{\geq}{\geqslant}

\def\epsilon{\varepsilon}
\def\phi {\varphi}
\def\supp{\text{supp}}

\begin{document}

\maketitle

\begin{abstract} In this work, we investigate the inverse problem of determining a quasilinear term appearing in a nonlinear elliptic equation from the measurement of the conormal derivative on the boundary. This problem arises in several practical applications, e.g., heat conduction. We derive novel H\"older stability estimates for both multi- and one-dimensional cases: in the multi-dimensional case, the stability estimates are stated with one single boundary measurement, whereas in the one-dimensional case, due to dimensionality limitation, the stability results are stated for the Dirichlet boundary condition varying in a space of dimension one. We derive these estimates using different properties of solution representations. We complement the theoretical results with numerical reconstructions of the quasilinear term, which illustrate the stable recovery of the quasilinear term in the presence of data noise.
\end{abstract}

\maketitle
\section{Introduction}
Let $\Omega \subset\R^n$ for $n \geq 1$  be an open  bounded domain with a ${C}^{2+\alpha}$ boundary $\partial\Omega$ for $\alpha \in (0,1)$, and $A(x)  = (a_{ij}(x))_{1\leq i,j \leq n}$, with $a_{ij}\in {C}^{1+\alpha}(\overline{\Omega}) $, be symmetric and elliptic:
\begin{equation}
    \label{eqn:ellipticity}
    a_{ij}(x) = a_{ji}(x), \, \xi^T A(x) \xi \geq c|\xi|^2, \quad x \in \overline{\Omega}, \,\xi \in \R^n.
\end{equation}
Consider the following quasilinear elliptic boundary value problem
\begin{equation}
    \label{eqn:ellptic problem}
    \left\{\begin{aligned}
        \nabla \cdot (\gamma (u)A(x)  \nabla u) &= 0,\quad \text{in } \Omega,\\
        u &= g,\quad \text{on } \partial\Omega.
    \end{aligned}\right.
\end{equation}
In this work we investigate the inverse problem of determining the quasilinear term $\gamma$ from the knowledge of the conormal derivative  of the solution $u$ of problem  \eqref{eqn:ellptic problem} restricted to a portion of the boundary $\partial\Omega$ for some suitable choice of $g$. For $n\geq2$, we investigate the problem with one single boundary measurement, i.e. the determination of $\gamma$ from the  measurement associated with one single Dirichlet excitation $g$. For $n=1$, due to the dimensionality limitation, it is no longer possible to have unique determination with one single boundary measurement, instead we consider Dirichlet excitations $g$ lying in a space of dimension one.

The concerned inverse problem arises in different physical and industrial applications, e.g., heat conduction, and describes the identification of a temperature dependent thermal conductivity $\gamma$ in a stationary diffusion process modeled by problem \eqref{eqn:ellptic problem} \cite{Al,BBC,Ca,CD1,Schuster}.
It has found applications in the manufacturing process of high-strength low-alloy water-cooled heavy plates, where the material properties are key targets \cite{Schuster}. The various macroscopic material properties are achieved by specific microstructures, and water cooling processes are an important instrument in adjusting
the microstructure \cite{SKFS}. The cooling processes however generate residual stresses that may adversely affect the flatness of the products, and the targeted
control of cooling processes on an industrial scale is imperative to meet the requirements on material properties
and flatness of the plates \cite{serizawa2015plate}.

In this work, we establish new conditional stability estimates for the inverse problem. In the two- and higher-dimensional case, we prove a H\"{o}lder stability estimate of recovering the quasilinear term  $\gamma$ from the Neumann data corresponding to one single Dirichlet excitation; see Theorem \ref{thm:main theoretical result} for the precise statement.  In the one-dimensional case, we prove a conditional stability result for the Neumann data corresponding to the Dirichlet excitation in a one-dimensional space for the quasilinear term $\gamma \in W^{1,p}((0,R))$, with $p>1$, and prove a uniqueness result for $\gamma \in W^{1,1}((0,R))$. See Theorems \ref{thm:1D-W1,1} and \ref{thm:1D-stability} for the precise statements. Moreover, we show that the stability in the 1D case is optimal in the sense that the H\"{o}lder exponent cannot be improved. Numerically we present regularized reconstructions for both one- and two-dimensional cases, and obtain accurate reconstructions  for noisy data.

The inverse problem of identifying  a nonlinear law described by quasilinear terms has a long history, and it has received much attention in the mathematical community in the last few decades. The first results related to this problem are stated in terms of unique determination of $\kappa(u)$ in the elliptic problem $-\nabla\cdot(\kappa(u)\nabla u)=f$ from the overspecified boundary data in the pioneering works of \cite{Ca,PR} (see also \cite{CD1,CD2} for the parabolic counterpart). More recently, there have been multiple works investigating this inverse problem either in the context of the minimal surface equation \cite{CLLO,Nu1,Nu2} or as an extension of the Calder\'on problem \cite{CFKKU,HS,KS,KKU,MU,Su,SuU}. These interesting works address the unique determination of quasilinear terms from the Dirichlet-to-Neumann map. So far in the literature,  there are only a few works that are devoted to the stable determination and the numerical reconstruction of quasilinear terms appearing in an elliptic equation from boundary measurements. Very recently, several researchers \cite{Cho,kian2023lipschitz,Ki24} proved stability estimates for the determination of a general class of quasilinear terms from boundary measurements associated with an infinite number of Dirichlet excitations (i.e., Dirichlet-to-Neumann type data); see also \cite{JK} for the parabolic counterpart. In contrast, the stability estimate in the multi-dimensional case in Theorem \ref{thm:main theoretical result} is stated for one single boundary measurement. To the best of our knowledge, for the one-dimensional case ($n=1$), the  determination of quasilinear terms from boundary measurements of solutions of elliptic equations has not been investigated in the literature.

In terms of numerical reconstruction of the quasilinear term, Kugler \cite{Ku} proposed a reconstruction method for recovering $\kappa(u)$ in $-\nabla\cdot(\kappa(u)\nabla u)=f$ from the overspecified boundary data based on standard Tikhonov regularization, and derived convergence rates for the regularized approximation using a problem-adapted adjoint under weaker conditions than the general theory \cite{EnglHankeNeubauer:1996}. Later Egger et al \cite{EPS4} proposed the Hilbert scale regularization for recovering $\kappa(u)$, based on a reformulation the inverse problem as a linear using the transform due to \cite{Ca} and an explicit characterization of mapping properties of the forward operator. More recently,  Kian et al \cite{KLWZ}  proved the stable determination and the numerical reconstruction of one class of semilinear terms in a semilinear elliptic equation from the boundary measurements. In this work we employ Tikhonov regularization \cite{EnglHankeNeubauer:1996,ItoJin:2015} for the numerical reconstruction of the quasilinear term.

The rest of the paper is organized as follows. In Section \ref{sec:main}, we present the mathematical formulation of the inverse problem, and state the main theoretical results. Then in Sections \ref{sec:proof-nd} and \ref{sec:proof-1d}, we present the proofs of the results for $n\geq2$ and $n=1$, respectively. In Section \ref{sec:numer} we present numerical results obtained by means of the standard regularized reconstruction. In the appendix we collect several technical results that are needed in the analysis of the one-dimensional case. Throughout we use standard notations for Sobolev spaces and continuously differentiable functions. The notation $c$ denotes a generic constant which may change from one line to the other, but it is always independent of the quasilinear term $\gamma$.

\section{Main results and discussions}\label{sec:main}
In this section we describe  preliminary tools and state the main results of this work. We first discuss the case $n \geq 2$. Consider problem \eqref{eqn:ellptic problem} with $\gamma \in {C}^2(\mathbb R;(0,+\infty))$ and  Dirichlet boundary data $ g \in {C}^{2+\alpha}(\partial\Omega)$, $\alpha\in(0,1)$,
with $\text{supp}(g) \subset S$, where $S\subset \partial\Omega$ is an arbitrary open and non-empty set. By \cite[p. 297, Theorem 8.2]{ladyzhenskaia1968linear} and \cite[Theorem 10.7]{gilbarg1977elliptic}, under the above assumptions,  problem \eqref{eqn:ellptic problem} has a unique solution $u \in {C}^{2+\alpha}(\overline{\Omega})$. Recall the tangential gradient $\nabla_\tau f$ for functions $f \in {C}^1(\partial\Omega)$ defined by
\begin{equation}
\label{eqn: defn of gradient for boundary fct.}
    \nabla_\tau f(x) = \nabla \tilde{f}(x) - (\partial_{\nu}\tilde{f}(x))\nu(x), \quad x \in \partial\Omega,
\end{equation}
where $\tilde{f}$ is in ${C}^1(\overline{\Omega})$ such that $\tilde{f}\mid_{\partial\Omega} = f$, and $\nu(x)$ is the unit outward normal vector to the boundary $\partial\Omega$ at $x \in \partial\Omega$. It is known that $\nabla_\tau f $ is independent of the choice of $\tilde{f}$. Then we can state the main result for the case $n \geq 2$, i.e., a H\"older stability result for recovering the quasilinear term $\gamma$.
\begin{theorem}
\label{thm:main theoretical result}\
    Let $n \geq 2$ and $R>0$. Suppose that $A$ satisfies \eqref{eqn:ellipticity} and  $\gamma_k \in {C}^2(\R)$ for $k=1,2$ satisfy the following: there exists constants $M>m>0$ such that
    \begin{equation}
    \label{eqn:cond on gamma}
    \|\gamma_k\|_{W^{1,\infty}((0,R))} \leq M \quad \mbox{and}\quad  \inf_{s \in [0,R] } \gamma_k(s) \geq m.
\end{equation}
Let the Dirichlet boundary data $g \in {C}^{2+\alpha}(\partial\Omega)$ be such that $\supp(g) \subset S$ with $S$ an arbitrary open and non-empty subset of $\partial\Omega$ and that $g(\partial\Omega) = [0,R]$. Additionally, assume that there exists $\delta \in [0,R/4)$ and $c_\delta>0$ such that for all $s \in [\delta,R-\delta]$, there exists $x_s\in S$ with
\begin{equation}
\label{eqn: cond on boundary fct. 1}
 g(x_s) = s \quad \mbox{and}\quad |\nabla_\tau g(x_s) | \geq c_\delta.
\end{equation}
Let $u_k\in {C}^{2+\alpha}(\overline{\Omega})$ be the solution to problem $\eqref{eqn:ellptic problem}$ with $\gamma=\gamma_k$. Then there exists a constant $c = c(g,\delta,m,M,R,\Omega,A)>0$ such that
\begin{equation}
\label{eqn: main theoretical result}
    \sup_{s\in [\delta,R-\delta]} |\gamma_1(s)-\gamma_2(s)| \leq c \|\gamma_1(u_1)\partial_{\nu_A} u_1 - \gamma_2(u_2)\partial_{\nu_A} u_2)\|^{\frac{1}{2+n}}_{L^1(S)}.
\end{equation}
\end{theorem}

To the best of our knowledge,  Theorem \ref{thm:main theoretical result} gives the first result of stable determination of a quasilinear term in a nonlinear elliptic equation from one single boundary measurement. Moreover, the single measurement can be located on an arbitrary open and nonempty portion $S$ of the boundary $\partial\Omega$ and the elliptic operator can have variable coefficients $A$. The only comparable results \cite{Cho,kian2023lipschitz,Ki24} provide similar stability estimates for a more general class of quasilinear terms, using an infinite number of measurements (i.e., Dirichlet-to-Neumann map). The use of one single measurement in the stability estimate \eqref{eqn: main theoretical result} makes the result especially suitable for the numerical reconstruction of the quasilinear term, when compared with the existing works that are using the Dirichlet-to-Neumann map. Indeed, using the stability estimate \eqref{eqn: main theoretical result}, we shall develop in Section \ref{sec:numer} a numerical algorithm for recovering the quasilinear term $\gamma$.

Next we consider the case $n=1$ and without loss of generality, we take $\Omega = (0,1)$. Fix $R>0$. Then problem \eqref{eqn:ellptic problem} takes the form
\begin{equation}
\label{eqn: bvp for 1d case}
    \left\{\begin{aligned}
                (\gamma (u)a(x)  u')' &= 0, \quad \mbox{in }(0,1),\\
        u(0) = 0, \, u(1) &= \lambda,
    \end{aligned}\right.
\end{equation}
where we take $\lambda \in [0,R]$, $a \in {C}^1([0,1])$ and $\gamma \in W^{1,p}((0,R);(0,+\infty))$, $p \in [1,\infty]$. We also assume that
\begin{equation}
\label{eqn: cond for 1d case looserW1p}
    \inf_{x\in[0,1]}a(x) > 0\quad\mbox{and} \quad \inf_{s\in[0, R]} \gamma(s) >0.
 \end{equation}
We shall show in Proposition \ref{prop:existence in 1d} that under condition \eqref{eqn: cond for 1d case looserW1p}, problem \eqref{eqn: bvp for 1d case} has a unique solution $u_\lambda$ lying in the Sobolev space $W^{2,p}((0,1);(0,R))$, which allows stating the following results in the case $n=1$. When $p=1$, we have a uniqueness result for recovering the quasilinear term $\gamma$:
\begin{theorem}
\label{thm:1D-W1,1}
    Let $R>0$ and suppose that $a \in {C}^1([0,1])$ and $\gamma_k\in W^{1,1}((0,R);(0,+\infty))$, $k=1,2,$  satisfy \eqref{eqn: cond for 1d case looserW1p}. Let $u_{k,\lambda} \in W^{2,1}((0,1);(0,R))$ be the solution to problem \eqref{eqn: bvp for 1d case} with $\gamma = \gamma_k$. Then, the condition
    $$\gamma_1(u_{1,\lambda}(1))a(1)u_{1,\lambda}'(1)=\gamma_2(u_{2,\lambda}(1))a(1)u_{2,\lambda}'(1), \quad \lambda \in [0,R],$$
    implies that $\gamma_1 = \gamma_2$ on $[0,R]$.
\end{theorem}

Meanwhile, for $p \in (1,\infty]$, we can prove a conditional H\"{o}lder stability result for the one-dimensional inverse problem. Moreover, in Remark \ref{rmk:optimality} below, we show that the H\"{o}lder exponent $\frac{p-1}{2p-1}$ in the stability estimate cannot be further improved.
\begin{theorem}
\label{thm:1D-stability}
    Let $R>0$, $p\in (1,\infty]$ and suppose that  $a \in {C}^1([0,1])$ and $\gamma_k\in W^{1,p}((0,R);(0,+\infty))$, $k=1,2,$ satisfy \eqref{eqn: cond for 1d case looserW1p}. Moreover, suppose that there exists $M>0$ such that
    \begin{equation}
              \label{eqn:cond on gamma for 1d W1p} M\geq\max(\|\gamma_1\|_{W^{1,p}((0,R))},\|\gamma_2\|_{W^{1,p}((0,R))}).
    \end{equation}
      Let $u_{k,\lambda} \in W^{2,p}((0,1);(0,R))$ be the solution to problem \eqref{eqn: bvp for 1d case} with $\gamma = \gamma_k$.  Then, there exists $c=c(a,p,M,R)>0$ such that
    \begin{equation}
    \label{eqn: stability for looser case}
    \|\gamma_1-\gamma_2\|_{L^\infty((0,R))} \leq c \sup_{\lambda\in[0,R]}|\gamma_1(u_{1,\lambda}(1))a(1)u_{1,\lambda}'(1)-\gamma_2(u_{2,\lambda}(1))a(1)u_{2,\lambda}'(1)|^{\frac {p-1}{2p-1}},
    \end{equation}
where $\frac {p-1}{2p-1}=\frac{1}{2}$ for $p=\infty$.
\end{theorem}

In contrast to the higher-dimensional case, Theorems \ref{thm:1D-W1,1} and \ref{thm:1D-stability} seem to be the first results dealing with the unique and stable determination of a quasilinear term from boundary measurements in dimension one. Due to the dimensionality limitation, it is impossible to uniquely determine a quasilinear term from one single boundary measurement comparable to that of Theorem \ref{thm:main theoretical result}. Therefore, in contrast to Theorem \ref{thm:main theoretical result}, we consider measurements associated with Dirichlet excitations lying in a space of dimension one. We obtain both sharp uniqueness and stability estimate stated for a general class of quasilinear terms where the regularity assumptions is drastically relaxed compared with the results in higher dimensions. Similar to Theorem \ref{thm:main theoretical result}, the flexibility of the results and the H\"older stability estimate in Theorem \ref{thm:1D-stability}, will be exploited in Section \ref{sec:numer} to derive the reconstruction algorithm.
Note that the proof of Theorem  \ref{thm:1D-stability} is based on the resolution of a variational problem of finding a minimizer of a functional subject to constraints. This approach allows us to derive an optimal H\"older exponent for the stability estimate \eqref{eqn:cond on gamma for 1d W1p} (see Remark \ref{rmk:optimality} for the proof of the optimality of the exponent $\frac{p-1}{2p-1}$, $p\in(1,+\infty]$).

Now we briefly comment on the overall proof strategies for the main theorems. The analysis for both the multi- and one-dimensional cases utilizes a key representation property of solutions of problem \eqref{eqn:ellptic problem} based on applying a suitable change of variables described in Section \ref{sec:proof-nd}. This representation is then combined with several properties of solutions of  the elliptic equations including regularity properties and the maximum principle. In contrast to most of the existing literature devoted to determination of quasilinear terms for an elliptic equation, see e.g. \cite{CFKKU,Cho,HS,kian2023lipschitz,Ki24,KKU,MU,Su,SuU}, we do not use the linearization technique, which necessarily involves infinite number of measurements. In the one-dimensional case, we first prove that up to a multiplicative constant, the antiderivative of $\gamma$ is equal to the Neumann data, cf. Proposition \ref{prop: main equality}. This key relation enables proving the uniqueness result in Theorem \ref{thm:1D-W1,1} directly. For the stability result in Theorem \ref{thm:1D-stability}, we construct a suitable variational problem and, under much relaxed regularity assumptions, employ some arguments used for the multidimensional problem with pointwise estimate for $\gamma$, leading
to the stability estimate in Theorem \ref{thm:1D-stability}.

\section{Proof of Theorem \ref{thm:main theoretical result}}\label{sec:proof-nd}
This section is devoted to the proof of  Theorem \ref{thm:main theoretical result}.

\begin{proof}[Proof of Theorem \ref{thm:main theoretical result}]
First, let
\begin{equation}
    \Gamma_k(s) = \int_0^s \gamma_k(t) \d t,
\end{equation}
and denote by $\mathcal{L}$ the elliptic operator defined by $\mathcal{L}v = \nabla \cdot(A \nabla v) $. Then in the domain $\Omega$, we have
$$\mathcal{L}(\Gamma_k(u_k))= \nabla \cdot (\gamma_k(u_k)A(x)\nabla u_k)=0.$$
Hence, by integration by parts,
$$
\begin{aligned}
    0=-&\int_\Omega (\mathcal{L}\Gamma_1(u_1)-\mathcal{L}\Gamma_2(u_2))(\Gamma_1(u_1)-\Gamma_2(u_2)) \d x \\
    =-&\int_{\partial\Omega} (\gamma_1(u_1)\partial_{\nu_A}u_1 - \gamma_2(u_2)\partial_{\nu_A}u_2)(\Gamma_1(u_1)-\Gamma_2(u_2)) \d\sigma(x) \\
    &+ \int_\Omega A\nabla(\Gamma_1(u_1)-\Gamma_2(u_2)) \cdot \nabla(\Gamma_1(u_1)-\Gamma_2(u_2))\d x.
\end{aligned}
$$ This implies the following identity
\begin{align}\label{t1d}
   &\int_\Omega A\nabla(\Gamma_1(u_1)-\Gamma_2(u_2)) \cdot \nabla(\Gamma_1(u_1)-\Gamma_2(u_2))\d x \\=&\int_{\partial\Omega} (\gamma_1(u_1)\partial_{\nu_A}u_1 - \gamma_2(u_2)\partial_{\nu_A}u_2)(\Gamma_1(u_1)-\Gamma_2(u_2)) \d\sigma(x).\nonumber
\end{align}
Note that $u_k=g$ on $\partial\Omega$, and $g(\partial\Omega) = [0,R]$. Moreover, since $\supp(g) \subset S$ and  $\Gamma_k(0) = 0$, $k=1,2$,  we deduce that $\Gamma_k(u_k) = 0$ on $\partial\Omega \setminus S$. By combining these facts with \eqref{eqn:cond on gamma} and \eqref{t1d}, we derive
\begin{align*}
&\left\lvert \int_\Omega A \nabla (\Gamma_1(u_1)- \Gamma_2(u_2))\cdot \nabla (\Gamma_1(u_1)-\Gamma_2(u_2)) \d x \right\lvert\\
=&\left\lvert\int_{S} (\gamma_1(u_1)\partial_{\nu_A}u_1 - \gamma_2(u_2)\partial_{\nu_A}u_2)(\Gamma_1(g(x)))-\Gamma_2(g(x))) \d\sigma(x)\right\lvert\\
\leq&2\left(\max_{j=1,2}\sup_{y\in S}|\Gamma_j(g(y))|\right)\int_S  |\gamma_1(u_1)\partial_{\nu_A}u_1 - \gamma_2(u_2)\partial_{\nu_A}u_2| \d\sigma(x)\\
\leq &2RM\int_S  |\gamma_1(u_1)\partial_{\nu_A}u_1 - \gamma_2(u_2)\partial_{\nu_A}u_2| \d\sigma(x).
\end{align*}
Moreover, by applying condition \eqref{eqn:ellipticity}, we conclude that there exists $c=c(M,R,\Omega,A)>0$ such that
\begin{equation}
\label{eqn: main thm proof est 1}
    \int_\Omega| \nabla (\Gamma_1(u_1(x))-\Gamma_2(u_2(x)))|^2 \d x \leq c\| \gamma_1(u_1)\partial_{\nu_A}u_1 - \gamma_2(u_2)\partial_{\nu_A}u_2\|_{L^1(S)}.
\end{equation}
Next, we derive  a lower bound for $|\nabla (\Gamma_1(u_1(x))-\Gamma_2(u_2(x)))|^2$. To this end, let $$\psi(x) := \nabla (\Gamma_1(u_1(x)) - \Gamma_2(u_2(x))) = \gamma_1(u_1(x))\nabla u_1(x) - \gamma_2(u_2(x)) \nabla u_2(x) . $$
Then, by taking $\delta$, $s$ and $x_s$ as in \eqref{eqn: cond on boundary fct. 1} and $\varepsilon >0$ to be arbitrary, and applying the mean value theorem to $\psi$, we  obtain
\begin{equation}
\label{eqn:est for psi}
    |\psi(x_s)-\psi(x)| \leq |\psi'(x_0)|\varepsilon, \quad x,x_0\in B(x_s,\varepsilon)\cap \Omega,
\end{equation}
where $B(x_s,\varepsilon)$ denotes the ball centered as $x_s$ with radius $\varepsilon$. By regarding $u_k$ as the solution to the linear boundary value problem:
$$
    \left\{\begin{aligned}
        \nabla \cdot (\gamma (u_k)A(x)  \nabla u) &= 0,\quad \text{in } \Omega,\\
        u &= g,\quad \text{on } \partial\Omega,
    \end{aligned}\right.
$$
we may apply the maximum principle for elliptic equations \cite[Chapter 6.4, Theorem 3]{evans2022partial} and deduce $u_k(\overline{\Omega}) \subset [0,R]$, $k=1,2$. Then, by  \cite[p. 270, Theorem 4.1]{ladyzhenskaia1968linear}
and \cite[p. 284, Theorem 6.5]{ladyzhenskaia1968linear},
we obtain further that $\| u_k\|_{C^2(\overline{\Omega})}\leq c$, $k=1,2$, for $c = c(g,m,M, R,\Omega,A)>0$. Thus the following estimate holds
$$
    \|\gamma(u_k) \nabla u_k\|_{C^1(\overline{\Omega})} \leq 3Mc,\quad k=1,2.
$$
Substituting the above estimate in \eqref{eqn:est for psi} leads to
$$
    |\psi(x_s)|^2 \leq (|\psi(x_s)-\psi(x)| + |\psi(x)|)^2 \leq 2(|\psi(x)|^2 + c\varepsilon^2), \quad x \in B(x_s,\varepsilon)\cap \Omega,
$$
which yields
$$
c|\psi(x)|^2 \geq |\psi(x_s)|^2 - c\varepsilon^2, \quad x \in B(x_s,\varepsilon)\cap \Omega.
$$
Integrating over $B(x_s,\varepsilon)\cap\Omega$ yields
\begin{equation}
\label{eqn: main thm proof est 2}
    c\int_\Omega |\psi(x)|^2 \d x  \geq c\int_{B(x_s,\varepsilon) \cap \Omega} |\psi(x)|^2 \d x \geq |B(x_s,\varepsilon) \cap \Omega||(\psi(x_s)|^2 - c\varepsilon^{2}).
\end{equation}
Moreover, since for $x_s\in\partial\Omega$,
we can find $\epsilon_0>0$ depending only on $\Omega$, such that for $\epsilon\in(0,\epsilon_0)$, there exists $c\in (0,1)$ depending only on $\Omega$ such that $|B(x_s,\varepsilon) \cap \Omega|=c\epsilon^n$.
Combining this with \eqref{eqn: main thm proof est 1} and \eqref{eqn: main thm proof est 2} yields
\begin{equation}
\label{eqn: main thm proof est 3}
    |\psi(x_s)|^2 \leq c(\varepsilon^{-n} \| \gamma_1(u_1)\partial_{\nu_A}u_1 - \gamma_2(u_2)\partial_{\nu_A}u_2\|_{L^1(S)} + \varepsilon^2),\quad \epsilon\in(0,\epsilon_0).
\end{equation}
Using the definition \eqref{eqn: defn of gradient for boundary fct.} of the tangential gradient $\nabla_\tau g$  and the boundary condition $u_k |_{\partial\Omega} = g$, we have
$$\nabla u_k = \nabla_\tau g + (\partial_\nu u_k) \nu,\quad \mbox{on }\partial\Omega.$$
Moreover, by assumption \eqref{eqn: cond on boundary fct. 1}, we have
$$u_k(x_s) = g(x_s) = s\quad\mbox{and}\quad |\nabla_\tau g(x_s)|\geq c_\delta.$$
Then the orthogonality relation $\nabla_\tau g \cdot \nu = 0$ implies
$$
\begin{aligned}
|\psi(x_s)|^2 &= |\gamma_1(s)(\nabla_\tau g(x_s) + \partial_{\nu}u_1(x_s)\nu(x_s)) - \gamma_2(s)(\nabla_\tau g(x_s) + \partial_{\nu}u_2(x_s)\nu(x_s))|^2 \\
&= |(\gamma_1(s) - \gamma_2(s))\nabla_\tau g(x_s) |^2 + |\gamma_1(s)\partial_\nu u_1(x_s) \nu(x_s) - \gamma_2(s) \partial_\nu u_2(x_s) \nu(x_s)|^2 \\
& \geq |\gamma_1(s)-\gamma_2(s)|^2 |\nabla_\tau g(x_s)|^2 \geq c_\delta^2|\gamma_1(s)-\gamma_2(s)|^2 .
\end{aligned}
$$
By combining this with the estimate \eqref{eqn: main thm proof est 3} and taking the supremum over all $s$, we get for  $\epsilon\in(0,\epsilon_0)$,
$$
\sup_{s \in [\delta,R-\delta]}|\gamma_1(s)-\gamma_2(s)|^2 \leq c(\varepsilon^{-n} \| \gamma_1(u_1)\partial_{\nu_A}u_1 - \gamma_2(u_2)\partial_{\nu_A}u_2\|_{L^1(S)} + \varepsilon^2),
$$
with $c = c(g,\delta,m,M,R,\Omega,A) >0$. Now, if $\| \gamma_1(u_1) \partial_{\nu_A} u_1 -\gamma_2(u_2) \partial_{\nu_A}u_2\|_{L^1(S)}= 0$, we have $$\sup_{s \in [\delta,R-\delta]}|\gamma_1(s)-\gamma_2(s)|^2 \leq c\varepsilon^2,$$
and we obtain the estimate \eqref{eqn: main theoretical result} by sending $\varepsilon$ to zero. For $\| \gamma_1(u_1) \partial_{\nu_A} u_1 -\gamma_2(u_2) \partial_{\nu_A}u_2\|_{L^1(S)}\in(0, \epsilon_0^{2+n})$, we may choose
$$\varepsilon = \| \gamma_1(u_1) \partial_{\nu_A} u_1 -\gamma_2(u_2) \partial_{\nu_A}u_2\|^{\frac{1}{2+n}}_{L^1(S)},$$
which yields the desired estimate \eqref{eqn: main theoretical result}. Finally, for
$\| \gamma_1(u_1) \partial_{\nu_A} u_1 -\gamma_2(u_2) \partial_{\nu_A}u_2\|_{L^1(S)}\geq \epsilon_0^{2+n}$, we can estimate directly
$$\sup_{s \in [\delta,R-\delta]}|\gamma_1(s)-\gamma_2(s)|\leq 2M\leq \frac{2M}{\epsilon_0}\| \gamma_1(u_1) \partial_{\nu_A} u_1 -\gamma_2(u_2) \partial_{\nu_A}u_2\|^{\frac{1}{2+n}}_{L^1(S)},$$
which again implies \eqref{eqn: main theoretical result}. This concludes the proof of the theorem.
\end{proof}

\section{Proofs of Theorems \ref{thm:1D-W1,1} and \ref{thm:1D-stability}}\label{sec:proof-1d}

In this section we present the proofs of Theorems \ref{thm:1D-W1,1} and \ref{thm:1D-stability}.
First we give the unique existence of a solution $u_\lambda\in W^{2,p}((0,1);[0,R])$ to problem \eqref{eqn: bvp for 1d case}.
\begin{proposition}
\label{prop:existence in 1d}
   Let $R>0$ and fix $\lambda \in [0,R]$, and suppose that $a \in {C}^1([0,1])$  and $\gamma \in W^{1,p}((0,R);(0,+\infty))$ for $p \in [1,\infty]$ satisfy the assumption
   \eqref{eqn: cond for 1d case looserW1p}. Then problem \eqref{eqn: bvp for 1d case} has a unique solution $u_\lambda \in W^{2,p}((0,1);[0,R])$.
\end{proposition}
\begin{proof}
First, note that any solution $u_\lambda$ must satisfy
$$
\gamma(u_\lambda(x))a(x) u_\lambda'(x) = c_1, \quad x \in [0,1],
$$
 where $c_1 \in \R$ is to be determined. By the monotonicity of the function $u_\lambda$, we have $0\leq u_\lambda\leq R$ and, by separation of variables, we derive that $u_\lambda$ must satisfy
$$
    \int_0^{u_\lambda(x)} \gamma(s) \d s = \int_0^x  \frac{c_1}{a(y)} \d y + c_2, \quad x \in [0,1],
$$
where $c_2 \in \R$ is again to be determined. We define a map $\Gamma \in {C}^1([0,R])$ by $\Gamma(s) = \int_0^{s} \gamma(s) \d s$. Indeed, for $p>1$, the continuous embedding of $W^{1,p}((0,R))$ into ${C}([0,R])$ holds  \cite[Theorem 8.8]{brezis2011functional}, and for $p=1$, functions in $W^{1,1}((0,R))$ are absolutely continuous \cite[Chapter 8, Remark 8]{brezis2011functional} and thus also lie in ${C}([0,R])$. Hence, $\gamma \in {C}([0,R])$ and $\Gamma \in C^1([0,R])$. Since $\Gamma'(\mu) = \gamma(\mu) > 0$ for all $\mu\in [0,R]$, by the inverse function theorem, $\Gamma$ is a $C^1$ diffeomorphism from $[0,R]$ to $[0,\Gamma (R)])$. Moreover, the condition $u_\lambda(0)=0$ yields $c_2 = 0$ and the condition $u_\lambda(1)=\lambda$ implies that
$$c_1=\left(\int_0^1\frac{1}{a(y)}  \d y \right)^{-1}\Gamma(\lambda)\geq0.$$
Thus, we have
$$0\leq\int_0^x\frac{c_1}{a(y)}  \d y\leq \int_0^1\frac{c_1}{a(y)}  \d y= \Gamma(\lambda)\leq \Gamma(R),\quad x\in[0,1]$$
and the unique solution $u_\lambda$ can be expressed as
    $$
    u_\lambda(x) = \Gamma^{-1} \left(\int_0^x\frac{c_1}{a(y)}  \d y \right),\quad x\in[0,1].
    $$
Clearly $u_\lambda\in C^1([0,1];(0,+\infty))$ and it remains to show that $u_\lambda'' \in L^p(0,1)$. By direct computation, we have
$$
u_{\lambda}''(x) = -\frac{c_1}{[a(x)]^2 \gamma(u_{\lambda}(x))} \left( a'(x) + \frac{c_1 \gamma'(u_{\lambda}(x))}{[\gamma(u_{\lambda}(x))]^2} \right)
$$
which implies that $u_\lambda''$ belongs to $L^p((0,1))$ by the assumptions on $a$ and $\gamma$ in \eqref{eqn: cond for 1d case looserW1p}. This concludes the proof of the proposition.
\end{proof}

The next result gives an important integral identity, which plays a crucial role in the analysis below.
\begin{proposition}\label{prop: main equality}
    Let $R>0$ and suppose $\lambda \in [0,R]$, and let $a$ and $\gamma$ be as in Proposition \ref{prop:existence in 1d}. Let $u_\lambda$ be the solution of problem \eqref{eqn: bvp for 1d case}. Define the Neumann data to  \eqref{eqn: bvp for 1d case} to be
\begin{equation}
\label{eqn:defintion of varphi}
    \varphi(\lambda) = \gamma(u_\lambda(1))a(1) u_\lambda'(1) .
    \end{equation}
    Then, for all $\lambda \in [0,R]$, there holds
  \begin{equation}
        \label{eqn: main equality in revised}
       \int_0^\lambda \gamma(s) \d s = c_a \varphi(\lambda),\quad \mbox{with }c_a = \int_0^1 \frac{1}{a(y)} \d y.
    \end{equation}
\end{proposition}
\begin{proof}
    First, fix $\lambda\in[0,R]$ and define
$\Gamma(s) = \int_0^s \gamma(\tau) \d \tau.$ Then we have
$$\frac{\d}{\d x}\Gamma(u_{\lambda}(x)) = \gamma(u_{\lambda}(x))u_{\lambda}'(x).
$$
By integrating equation \eqref{eqn: bvp for 1d case} and applying the fundamental theorem of calculus, we get
\begin{align*}0&=\int_1^x\frac{\d}{\d y}
    \left(a(y)\frac{\d}{\d y}(\Gamma(u_{\lambda}(y)))\right)\d y \\
    &=a(x)\frac{\d}{\d x}(\Gamma(u_{\lambda}(x)))-a(1)\frac{\d}{\d x}(\Gamma(u_{\lambda}(1)))\\
    &=a(x)\frac{\d}{\d x}(\Gamma(u_{\lambda}(x)))- \gamma(\lambda)a(1)u_{\lambda}'(1),\quad x\in[0,1].
\end{align*}
Then it follows that

    \begin{equation}
  \label{ts1}
    \frac{\d}{\d x}(\Gamma(u_{\lambda}(x))) = \frac{1}{a(x)}\varphi(\lambda),\quad x\in[0,1].
  \end{equation}
  Upon noting the condition $u_{\lambda}(0) =0$, we deduce $\Gamma(u_{\lambda}(0))=\Gamma(0) = 0$. Then by integrating the identity \eqref{ts1} with respect to $x$, we obtain
$$
        \Gamma(u_{\lambda}(1))=\int_0^1\frac{\d}{\d x}(\Gamma(u_{\lambda}(x)))\d x  = \varphi(\lambda) \int_0^1 \frac{1}{a(y)} \d y.
$$
Upon noting the condition $u_\lambda(1) = \lambda$ and also using the identity $\Gamma(\lambda)=\int_0^\lambda \gamma(s)\d s$, the desired identity \eqref{eqn: main equality in revised} follows directly.
\end{proof}

With Proposition \ref{prop: main equality}, we can now prove both Theorems \ref{thm:1D-W1,1} and \ref{thm:1D-stability}.

\begin{proof}[Proof of Theorem \ref{thm:1D-W1,1}]

Let \( \varphi = \varphi_1 - \varphi_2 \), where \( \varphi_k \) ($k=1,2$) is the Neumann data defined in \eqref{eqn:defintion of varphi} with \( \gamma = \gamma_k \), and also let \( \gamma = \gamma_1 - \gamma_2 \). Let $I=\{s\in[0,R]:\ |\gamma(s)|>0\}$. Since $\gamma \in W^{1,1}((0,R)) \subset C([0,R])$, either ${I}=\emptyset$ or ${I}$ contains a nonempty open interval. We assume for a contradiction that $I \neq \emptyset$. By assumption, there exists $s_1,s_2\in[0,R]$, $s_1<s_2$ such that $(s_1,s_2) \subset {I}$. By the intermediate value theorem we can assume without loss of generality that, for all $s\in(s_1,s_2)$, $\gamma(s)>0$ and by \eqref{eqn: main equality in revised}, we have
    $$
    \begin{aligned}
    0<\int_{s_1}^{s_2} \gamma(s) \d s &= \int_0^{s_2} (\gamma_1(s) - \gamma_2(s)) \d s - \int_0^{s_1} (\gamma_1(s) - \gamma_2(s)) \d s \\ &= c_a[\varphi_1(s_2) - \varphi_2(s_2) - (\varphi_1(s_1) - \varphi_2(s_1))] = 0.
    \end{aligned}
    $$
    This leads to a contradiction and we deduce that  ${I}=\emptyset$ which implies $\gamma_1 = \gamma_2$ on $[0,R]$.
\end{proof}

The proof of Theorem \ref{thm:1D-stability} is lengthy and thus several technical results will be deferred to Appendix \ref{sec:appendix}.

\begin{proof}[Proof of Theorem \ref{thm:1D-stability}]
Again let \( \varphi = \varphi_1 - \varphi_2 \), with \( \varphi_k \) ($k=1,2$) the Neumann data defined in \eqref{eqn:defintion of varphi} with \( \gamma = \gamma_k \), and also let \( \gamma = \gamma_1 - \gamma_2 \). Now we treat the two cases $1<p<\infty$ and $p=\infty$ separately. We will first treat the case $1<p<\infty$. The proof is lengthy and technical, and we split it into several steps.

\noindent\textbf{Step 1: Construct suitable variational problem.}
 Let $I=(s_1,s_2) \subset\{s \in [0,R]:|\gamma(s)|>0\}$ be any open subinterval of $(0,R)$ on which $\gamma$ is non-vanishing. By the intermediate value theorem we assume without loss of generality that $\gamma > 0$ on $I$. Then, by the assumption \eqref{eqn:cond on gamma for 1d W1p} and the triangle inequality, we have
\begin{equation}
\label{eqn: ts2 p1}
\quad \|\gamma'\|_{L^p(I)} \leq 2M,
\end{equation}
and by Proposition \ref{prop: main equality}, we have
\begin{equation}
    \label{eqn: ts2 p2}
\| \gamma \|_{L^1(I)}  = \int_I \gamma_1(s) - \gamma_2(s) \d s = c_a(\varphi(s_2) - \varphi(s_1)) \leq 2c_a \| \varphi\|_{L^\infty(I)}.
\end{equation}
It hence suffices to prove that for all possible choices of $I$, the following estimate holds
\begin{equation}
\label{eqn:ineq_p_finite}
\left( \frac{p-1}{2p-1} \right)^{p-1} \frac{\|\gamma\|_{L^\infty(I)}^{2p-1}}{\|\gamma\|_{L^1(I)}^{p-1}} \leq \|\gamma'\|_{L^p(I)}^p.
\end{equation}
This and the estimates \eqref{eqn: ts2 p1} and \eqref{eqn: ts2 p2} yield directly
$$
\|\gamma\|_{L^\infty((0,R))}^{2p-1} \leq 2^{2p-1} M^p c_a^{p-1} \left( \frac{2p-1}{p-1} \right)^{p-1}  \|\varphi\|_{L^\infty((0,R))}^{p-1},
$$
which directly implies \eqref{eqn: stability for looser case} for $1<p<\infty$. It suffices to consider the cases where $I=(s_1,s_2)$ satisfies one of the three conditions below:
 $$
\gamma(s_1)=\gamma(s_2) = 0; \quad \gamma(s_1)=0,\, s_2=R, ; \quad s_1 = 0, \, s_2 = R;
 $$
and the case $s_1=0$, $\gamma(s_2) = 0$ follows from symmetry. We will first prove \eqref{eqn:ineq_p_finite} for the cases $\gamma(s_1)=\gamma(s_2)=0$ and $\gamma(s_1)=0$, $s_2=R$ by using a suitable variational problem. The third case $s_1=0$, $s_2=R$ will be discussed below in Step 4. Specifically, consider the problem of minimizing \( \|f'\|_{L^p(I)}^p \) over the space \(W^{1,p}(I) \) subject to the constraints \( \|f\|_{L^\infty(I)} = H \), \( \|f\|_{L^1(I)} = S \), \( f \geq 0 \), and \( \min_{x \in \partial I} f = 0 \). Define
\begin{equation}
\label{eqn: defn of M^p_H,S}
\mathfrak{M}^p_{I,H,S} = \left\{ f \in W^{1,p}(I) \mid \|f\|_{L^\infty(I)} = H, \|f\|_{L^1(I)} = S, f \geq 0, \min_{x \in \partial I} f = 0 \right\},
\end{equation}
and the functional
\begin{equation}
\label{eqn:defn of F_I^p}
\mathcal{F}^p_I(f) = \|f'\|_{L^p(I)}^p.
\end{equation}
To prove the estimate \eqref{eqn:ineq_p_finite}, it suffices to show that the minimizer \( f_0 = \arg\min_{f \in \mathfrak{M}^p_{I,H,S}} \mathcal{F}^p_I(f) \) exists and satisfies
\[
\left( \frac{p-1}{2p-1} \right)^{p-1} \frac{H^{2p-1}}{S^{p-1}} \leq \mathcal{F}^p_I(f_0).
\]

\noindent\textbf{Step 2: The existence of a minimizer}. Let $\{f_n\}_{n=1}^\infty \subset \mathfrak{M}^p_{I,H,S}$ be a minimizing sequence for $\mathcal{F}^p_I$.  
Moreover, using that $\|f_n\|_{L^\infty(I)} = H$, we deduce that the sequence $\{\|f_n\|_{W^{1,p}(I)}\}_{n=1}^\infty$ must be uniformly bounded. Note that by Sobolev embedding theorem \cite[Theorem 8.8]{brezis2011functional}, $W^{1,p}(I)$ is compactly embedded into $C(\overline{I})$. Thus by the Banach-Alaoglu theorem \cite[Theorem 3.16]{brezis2011functional}, there exists a weakly convergent subsequence of $\{f_n\}_{n=1}^\infty$ (still denoted by $\{f_n\}_{n=1}^\infty$) and $f_0 \in W^{1,p}(I)$ in the sense of $W^{1,p}(I)$ such that weakly in $W^{1,p}(I)$ to $f_0$ and $\{f_n\}_{n=1}^\infty$ converges strongly in  ${C}(\overline{I})$. Now, \cite[Chapter 8.2, Theorem 1]{evans2022partial}, the functional $\mathcal{F}^p_I$ is weakly lower semicontinuous, which implies
$$
\mathcal{F}_I^p(f_0) \leq \liminf_{n\to\infty} \mathcal{F}_I^p(f_n).
$$
Moreover, since $f_n \to f_0$ in ${C}(\overline{I})$, it follows directly that
$$
\|f_0\|_{L^\infty(I)}=H, \quad \|f_0\|_{L^1(I)}=S, \quad f_0 \geq 0, \quad \min_{x \in \partial I} f_0 = 0,
$$
and hence $f_0 \in \mathfrak{M}_{I,H,S}^p$. It follows that $f_0$ is a minimizer for $\mathcal{F}_I^p$ under the given constraints.

\noindent \textbf{Step 3: Estimate $\mathcal{F}_I^p(f_0)$}. First we find an explicit expression for the minimizer $f_0$. Suppose first that $f_0$ satisfies \(0< f_0 < H\) on an open interval \(I' \subset I\). Then, for any function \(f \in C_c^\infty(I')\) such that \(\int_{I'} f \, \d s = 0\), we have
$$
\left. \frac{\d}{\d \tau}\mathcal{F}^p_I(f_0 + \tau f)\right\lvert_{\tau = 0}   = 0.
$$
Now, by regarding $f_0$ as a distribution on $I'$, we derive
\[
\left\langle (p(f_0')^{p-1} )',f \right\rangle_{D'(I'),C^\infty_c(I')} = - \int_{I'} p(f_0')^{p-1} f' \d s = \left. \frac{\d}{\d \tau}\mathcal{F}^p_I(f_0 + \tau f)\right\lvert_{\tau = 0} = 0.
\]
Thus $(p(f_0')^{p-1} )'$ must be constant on $I'$. Hence, there exist constants $b$, $\tilde{b}$ and $c$ such that
\begin{equation}
\label{eqn: form of f0 on I'}
f_0(s)= c\left(s+b\right)^{\frac{p}{p-1}}+\tilde{b},\quad \mbox{on } I'.
\end{equation}
 Now, note that we can decompose the internal $I=(s_1,s_2)$ into $ I = \{f_0 = 0\} \cup \{0<f_0<H\}\cup \{f_0=H\} $. Since the set $\{0<f_0<H\}$ is a union of open intervals $I'$, $f_0$ must be defined piecewise by \eqref{eqn: form of f0 on I'}. Furthermore, we can derive an explicit form of $f_0$ and thus estimate $\mathcal{F}^p_I(f_0)$ from below. The details are technical and lengthy, and thus the derivation is deferred to the appendix, specifically in Propositions \ref{prop:Needyform1}$\sim$\ref{prop:minimal_F}. In Propositions \ref{prop:Needyform1} and \ref{prop:Needyform2}, the explicit form of $f_0$ under the respective additional restrictions $f(s_1) = f(s_2) = 0$ and $f(s_1) = 0$, $s_2 = R$ are given; and in Proposition \ref{prop:minimal_F}, the following estimate $$\mathcal F_I^p(f_0) \geq   \left(\frac{p-1}{2p-1}\right ) ^{p-1}\frac{ H^{2p-1}}{S^{p-1}}$$
 is established. This proves the claim \eqref{eqn:ineq_p_finite} in the cases $\gamma(s_1) = \gamma(s_2) = 0$ and $\gamma(s_1) = 0$, $s_2 = R$.

\noindent\textbf{Step 4: The remaining case that $s_1=0$, $s_2=R$.} It suffices to prove the estimate \eqref{eqn:ineq_p_finite} for $\gamma$ with \( \gamma > 0 \) on \( [0, R] \). Consider a vertical translation of $\gamma_1$, namely $\gamma_3 = \gamma_1 - h$ with $h = \min_{s \in [0,R]} \gamma(s)$. By applying \eqref{eqn: ts2 p2}, it follows that
$$
0 < h \leq R^{-1} \|\gamma\|_{L^1((0,R))} \leq 2c_a R^{-1} \|\varphi\|_{L^\infty((0,R))}.
$$
Letting $s_c =  \arg\min_{s \in [0,R]} \gamma(s)$, we see that $\gamma_3(s_c) = \gamma_2(s_c)$. We must also have $\gamma_3 > 0$, in view of the identity $\gamma_3 = \gamma_2 + (\gamma - h) >0$. Now, by the assertion in Step 3, we may take $\gamma = \gamma_3 - \gamma_2$ and $\varphi = \varphi_3 - \varphi_2$ (with $\varphi_3(s) = c_a^{-1}\int_0^s \gamma_3(t) \d t$) in \eqref{eqn:ineq_p_finite} with $I = (0,s_c)$ or $I = (s_c,R)$ and obtain
$$\|\gamma_3-\gamma_2\|_{L^\infty((0,R))}^{2p-1}\leq c\|\varphi_3-\varphi_2\|^{p-1}_{L^\infty((0,R))}\leq c(\|\varphi_1 - \varphi_2\|_{L^\infty((0,R))}+h)^{p-1}.$$
Then by the triangle inequality,
$$\|\gamma_1 - \gamma_2\|_{L^\infty((0,R))}\leq \|\gamma_3-\gamma_2\|_{L^\infty((0,R))}+h\leq c(\|\varphi_1 - \varphi_2\|_{L^\infty((0,R))} + h)^{\frac{p-1}{2p-1}} + h$$
Since $h \leq 2c_aR^{-1}\|\varphi\|_{L^\infty((0,R))}$, we get
$$
\|\gamma_1 - \gamma_2\|_{L^\infty((0,R))} \leq c\Big(\|\varphi_1 - \varphi_2\|_{L^\infty((0,R))}^\frac{p-1}{2p-1} + \|\varphi_1 - \varphi_2\|_{L^\infty((0,R))}\Big) \leq c \|\varphi_1 - \varphi_2\|_{L^\infty((0,R))}^\frac{p-1}{2p-1}.
$$
This concludes the proof of the claim \eqref{eqn:ineq_p_finite} and thus we have shown that the estimate \eqref{eqn: stability for looser case} holds for $1<p<\infty$.

\noindent\textbf{(Step 5): The case $p=\infty$.} Last, we prove the case $p=\infty$, which will be proved directly.
By Sobolev embedding, we have $\gamma_k \in W^{1,\infty}((0,R)) \subset {C}([0,R])$. Hence we can find $s_c \in [0,R]$ such that $\|\gamma\|_{L^\infty((0,R))} = |\gamma(s_c)|$. Without loss of generality, we suppose that $\gamma(s_c) \geq 0$. Let $0<\varepsilon < R/2$ be arbitrary and define $s_1 = \max(0,s_c-\varepsilon)$, $s_2 = \min(R,s_c+\varepsilon)$. Then, by applying \eqref{eqn: main equality in revised}, we have
    $$
    \begin{aligned}
    \int_{s_1}^{s_2} \gamma(s) \d s = \int_0^{s_2}\gamma(s) \d s - \int_0^{s_1} \gamma(s) \d s
    =c_a(\varphi(s_2) - \varphi(s_1))
    \leq 2c_a \sup_{\lambda \in [0,R]} |\varphi(\lambda)|.
    \end{aligned}
    $$
By the mean value theorem, we obtain that for any $s \in (s_1,s_2)$
    $$
    \gamma(s_c) \leq |\gamma(s_c) - \gamma(s)| + \gamma(s) \leq 2M\varepsilon + \gamma(s).
    $$
    Integrating from $s_1$ to $s_2$ and applying the previous estimate yields
    $$
    (s_2-s_1)\gamma(s_c) \leq 2M(s_2-s_1)\varepsilon + 2c_a \sup_{\lambda \in [0,R]} |\varphi(\lambda)|.
    $$
    Now, upon noting the condition $s_2-s_1 \geq \varepsilon$, we obtain
    $$
    \gamma(s_c) \leq c\left(\varepsilon^{-1}\sup_{\lambda \in [0,R]}|\varphi(\lambda)| + \varepsilon\right).
    $$
    If $\sup_{\lambda\in [0,R]}|\varphi(\lambda)|=0$, the estimate \eqref{eqn: stability for looser case} follows by sending $\varepsilon$ to zero. If $\sup_{\lambda\in [0,R]}|\varphi(\lambda)|^{\frac{1}{2}} \in  (0,R/2)$, the estimate \eqref{eqn: stability for looser case} follows by taking  $\varepsilon = \sup_{\lambda\in [0,R]}|\varphi(\lambda)|^{\frac{1}{2}}$. Lastly, for $\sup_{\lambda\in [0,R]}|\varphi(\lambda)|^{\frac{1}{2}} \geq R/2$, we have
$$
|\gamma(s_c)| = \gamma(s_c) \leq 2M\leq \frac{4M}{R}\sup_{\lambda\in [0,R]}|\varphi(\lambda)|^{\frac{1}{2}} \leq c\sup_{\lambda\in[0,R]}|\varphi(\lambda)|^{\frac{1}{2}}.
$$
This concludes the proof of the case $p=\infty$ and also the theorem.
\end{proof}

\begin{remark}\label{rmk:optimality}
The H\"older exponent $\frac{p-1}{2p-1}$ in Theorem \ref{thm:1D-stability} is optimal for all $p \in (1,\infty]$. To see this, fixing $a \equiv 1$ and $R = 1$, we claim that for all $q > \frac{p-1}{2p-1}$ $(\frac{p-1}{2p-1} = \frac{1}{2}$ in the case $p=\infty)$ and for all $c=c(M,p)>0$, there exist $\gamma_1,\gamma_2 \in W^{1,p}((0,1),(0,+\infty))$ satisfying \eqref{eqn: cond for 1d case looserW1p} and \eqref{eqn:cond on gamma for 1d W1p} for some fixed $M>0$ such that
\begin{equation}
\label{eqn:optimality}
    \|\gamma_1-\gamma_2\|_{L^\infty((0,1))} > c \sup_{\lambda \in [0,R]} |\varphi(\lambda)|^q,
\end{equation}
    with $\varphi(\lambda) = \varphi_1(\lambda)-\varphi_2(\lambda)$. We define $\gamma_1$ by
\[
\gamma_1(s) = \begin{cases}
1, & 0 \leq s \leq s_c, \\
1 + \frac{H}{1-s_c}(s - s_c) , & s_c < s \leq 1,
\end{cases}
\]
where $H$ and $s_c$ are to be determined, and take  $\gamma_2 \equiv 1$. Then there holds
$$
\gamma_1(s) - \gamma_2(s) = \begin{cases}
    0, & 0 \leq s \leq s_c, \\
     \frac{H}{1-s_c}(s - s_c), & s_c<s\leq 1,\\
\end{cases}
$$
and by applying \eqref{eqn: main equality in revised}, we get
$$
\varphi(\lambda) = \int_0^\lambda \gamma_1(s) -\gamma_2(s) \d s = \begin{cases}
    0, & 0\leq \lambda \leq s_c, \\
   \frac{H}{2( 1-s_c)}(\lambda-s_c)^2 & s_c<\lambda \leq 1.
\end{cases}
$$
Now, let $0<S<2^{-\frac{2p-1}{p}}M^{\frac{p-1}{p}}$ ($-\frac{2p-1}{p} = -2$ and $\frac{p-1}{p} = 1$ in the case $p=\infty$) be arbitrary, and set $H=(MS)^{\frac{p-1}{2p-1}}$ and $ s_c = 1 - \frac{2S}{H}$. Note that the conditions on $S$ imply that $0<s_c<1$. Direct computation gives
$$
\|\gamma_1'\|_{L^p((0,1))}^p = \frac{H^p}{(1-s_c)^{p-1}} = \frac{H^{2p-1}}{(2S)^{p-1}}=\left(\frac{M}{2}\right)^{p-1}
$$
for $1< p < \infty$, and $\|\gamma_1'\|_{L^\infty((0,1))} \leq \frac{M}{2}$ in the case $p=\infty$. Thus, $\gamma_1$ and $\gamma_2$ satisfy the condition \eqref{eqn:cond on gamma for 1d W1p}. Moreover, we have
\[
\|\gamma_1 - \gamma_2 \|_{L^\infty((0, 1))} = H =(MS)^{\frac{p-1}{2p-1}}, \quad \sup_{\lambda \in [0,R]} |\varphi(\lambda)|= S.
\]
Thus, we can choose $S =\min((M^{\frac{p-1}{2p-1}}c^{-1})^{\left(q-\frac{p-1}{2p-1}\right)^{-1}},2^{-\frac{2p-1}{p}}M^{\frac{p-1}{p}})$. With this choice, the functions $\gamma_1$ and $\gamma_2$ fulfill \eqref{eqn:optimality}.
\end{remark}

\section{Numerical experiments and discussions}\label{sec:numer}

In this section, we present numerical reconstructions to illustrate the stability of the inverse problem in both one- and two-dimensional cases.

\subsection{Numerical algorithm}
\label{subsec:reconstruction algorithm}
First we describe the numerical algorithm. The Neumann boundary data \(\gamma(u_\gamma) \partial_{\nu_A} u_\gamma\) is given by
$v = \gamma(g) \left. \left(A \nabla u_\gamma(g)\right) \right|_{\partial \Omega} \cdot \nu$,
where \(u_\gamma: g \mapsto u_\gamma\) is the solution operator for problem \eqref{eqn:ellptic problem}. We employ the standard Tikhonov regularization \cite{EnglHankeNeubauer:1996,ItoJin:2015} with an $H^1((0,R))$ semi-norm penalty for the  reconstruction. It amounts to minimizing the following regularized functional
\begin{equation}
    {J}_\beta(\gamma) =  {\frac{1}{2} \int_{S} (\gamma (u_\gamma) (\partial_{\nu_A}u_\gamma - v))^2\d\sigma}
+ {\frac{\beta}{2} |\gamma|_{H^1((0,R))}^2,   }
\label{eq:objective_functional}
\end{equation}
where $\beta>0$ is the regularization parameter balancing the data fidelity and regularization term. The regularization term $|\gamma|_{H^1((0,R))}$ is to stabilize the minimization problem so that the quasilinear term $\gamma$ has suitable regularity. The absence of the term leads to numerical instability, especially when the initial guess is far from the exact solution.

To minimize the functional ${J}_\beta(\gamma)$, we employ the gradient descent method. Below we derive a formula for the Gateaux derivative ${J}_0'(\gamma)[ h]=\frac{\d}{\d\tau}({J}_0(\gamma+\tau h))|_{\tau =0}$ of the data-fitting term ${J}_0$ using the standard adjoint technique. The next result shows that we can compute the gradient by solving  for $z$ a linear elliptic equation. The notation $(\cdot,\cdot)_{L^2(\Omega)}$ denotes the $L^2(\Omega)$ inner product, and similarly $(\cdot,\cdot)_{L^2(\partial\Omega)}$.
\begin{proposition}\label{prop:derivative}
For the functional ${J}_0(\gamma)$, the Gateaux derivative ${J}_0'(\gamma)$ is given by
\begin{equation*}
     {J}_0'(\gamma)[h]= \langle\nabla z\cdot A\nabla u_\gamma,h(u_\gamma)\rangle_{L^2(\Omega)},
\end{equation*}
where the function $z$  satisfies the adjoint problem
$$
    \left\{\begin{aligned}
         \nabla\cdot( A \nabla z) &=0, \quad\mbox{in }\Omega,\\
        z&=\chi_S(\gamma(u_\gamma)\partial_{\nu_A} u_\gamma-v) ,\quad \mbox{on }\partial\Omega.
    \end{aligned}\right.
$$
\end{proposition}
\begin{proof}
Note that the directional derivative $w = \left(\frac{\d }{\d \tau}u_{\gamma+\tau h}\right)\mid_{\tau = 0}$ satisfies
\begin{equation}
    \label{eqn:eqn-for-w}
    \left\{\begin{aligned}
        \nabla\cdot\big((\gamma'(u_\gamma) w + h(u_\gamma))A\nabla u_\gamma + \gamma(u_\gamma)A\nabla w)\big)&=0,\quad  \text{in } \Omega,\\
        w&=0,\quad \text{on } \partial\Omega.
    \end{aligned}\right.
\end{equation}
Then it follows directly from the product rule and the zero boundary condition $w=0$ on $\partial\Omega$ that
\begin{align}
    {J}_0'(\gamma)[h]&= \int_S(\gamma(u_\gamma)\partial_{\nu_A} u_\gamma-v)\big(\gamma'(u_\gamma)w+h(u_\gamma))\partial_{\nu_A} u_\gamma+\gamma(u_\gamma)\partial_{\nu_A} w\big)\d\sigma\nonumber\\
    & =\int_S(\gamma(u_\gamma)\partial_{\nu_A} u_\gamma-v)\big(h(u_\gamma)\partial_{\nu_A} u_\gamma +\gamma(u_\gamma)\partial_{\nu_A} w\big)\d\sigma\nonumber\\
    &=\langle  z,h(u_\gamma)\partial_{\nu_A}  u_\gamma\rangle _{L^2(\partial\Omega)}+  \langle z,\gamma(u_\gamma)\partial_{\nu_A} w\rangle_{L^2(\partial\Omega)}.\label{eqn:derivative-0}
\end{align}
Then by Green's identity and integration by parts (using $w=0$ on $\partial\Omega$ again), we derive
\begin{align*}
     \langle z,\gamma(u_\gamma)\partial_{\nu_A} w\rangle_{L^2(\partial\Omega)}
    &=\langle z,\nabla\cdot(\gamma(u_\gamma)A\nabla w)\rangle_{L^2(\Omega)}+\langle \nabla z,\  \gamma(u_\gamma)A\nabla w \rangle_{L^2(\Omega)}\\
   & = \langle z,\nabla\cdot(\gamma(u_\gamma)A\nabla w)\rangle_{L^2(\Omega)}- \langle \nabla\cdot(\gamma(u_\gamma)A \nabla z), w \rangle_{L^2(\Omega)} .
\end{align*}
Using the defining identity \eqref{eqn:eqn-for-w} for $w$ and integration by parts, we have
\begin{align*}
\langle z,\nabla\cdot(\gamma(u_\gamma)A\nabla w)\rangle_{L^2(\Omega)}&=-\langle z, \nabla\cdot\big(h(u_\gamma)A\nabla u_\gamma+\gamma'(u_\gamma)w A\nabla u_\gamma\big)\rangle_{L^2(\Omega)}\\
    &=-\langle z, \nabla\cdot\big(h(u_\gamma)A\nabla u_\gamma\big)\rangle_{L^2(\Omega)}+\langle\nabla z,\gamma'(u_\gamma)w A\nabla u_\gamma\rangle_{L^2(\Omega)}.
\end{align*}
Now using the elementary identity
$$\langle \nabla\cdot(\gamma(u_\gamma)A \nabla z), w \rangle_{L^2(\Omega)} =
\langle \gamma'(u_\gamma)\nabla u_\gamma \cdot A \nabla z, w \rangle_{L^2(\Omega)} +\langle \gamma (u_\gamma)\nabla\cdot(A \nabla z), w \rangle_{L^2(\Omega)}, $$
and combining the preceding estimates with the defining identity of $z$ yields
\begin{align*}
 \langle z,\gamma(u_\gamma)\partial_{\nu_A} w\rangle_{L^2(\partial\Omega)}
  & =-\langle z, \nabla\cdot\big(h(u_\gamma)A\nabla u_\gamma\big)\rangle_{L^2(\Omega)}-\langle\gamma(u_\gamma)\nabla\cdot(A\nabla z),w\rangle_{L^2(\Omega)}\\
   &=-\langle z, \nabla\cdot\big(h(u_\gamma)A\nabla u_\gamma\big)\rangle_{L^2(\Omega)}\\
   &=\langle\nabla z,h(u_\gamma)A\nabla u_\gamma\rangle_{L^2(\Omega)}-\langle  z,h(u_\gamma)\partial_{\nu_A}  u_\gamma\rangle _{L^2(\partial\Omega)}.
\end{align*}
Thus by combining the identity with \eqref{eqn:derivative-0}, we obtain the desired expression of the gradient ${J}_0'(\gamma)[h]$.
\end{proof}

Using Proposition \ref{prop:derivative}, one can use any gradient type method to minimize the objective functional ${J}_\beta(\gamma)$, e.g., gradient descent, stochastic gradient descent or Adam algorithm \cite{Kingma2014AdamAM}.

\subsection{Numerical results and discussions}

Now we numerically investigate the feasibility of reconstructing a quasilinear term. The regularizer $|\gamma|_{H^1(0,R)}^2$ is taken to be
$$|\gamma|_{H^1(0,R)}^2=\sum_{k=1}^{N-1} \int_{u_\gamma(x_{k})}^{u_\gamma(x_{k+1})}|\gamma'(s)|^2\d s,$$
where $\{x_k\}$ are mesh points on the boundary, and over each sub-interval, we employ linear interpolation for $\gamma$. The exact Neumann data $v$ is obtained by solving the direct problem \eqref{eqn:ellptic problem} on a fine  mesh, and the inverse problem is solved using a coarse mesh, in order to avoid inverse crime. The noisy  data $v^\delta$ is generated by adding Gaussian random noise multiplicatively into $v$.
The resulting optimization problem is solved using a gradient descent type method.

\begin{figure}[htb!]
    \begin{subfigure}[b]{\linewidth}
    \centering
    \begin{tabular}{c@{\hspace{2mm}}c@{\hspace{2mm}}c}
           \includegraphics[width=0.30\linewidth]{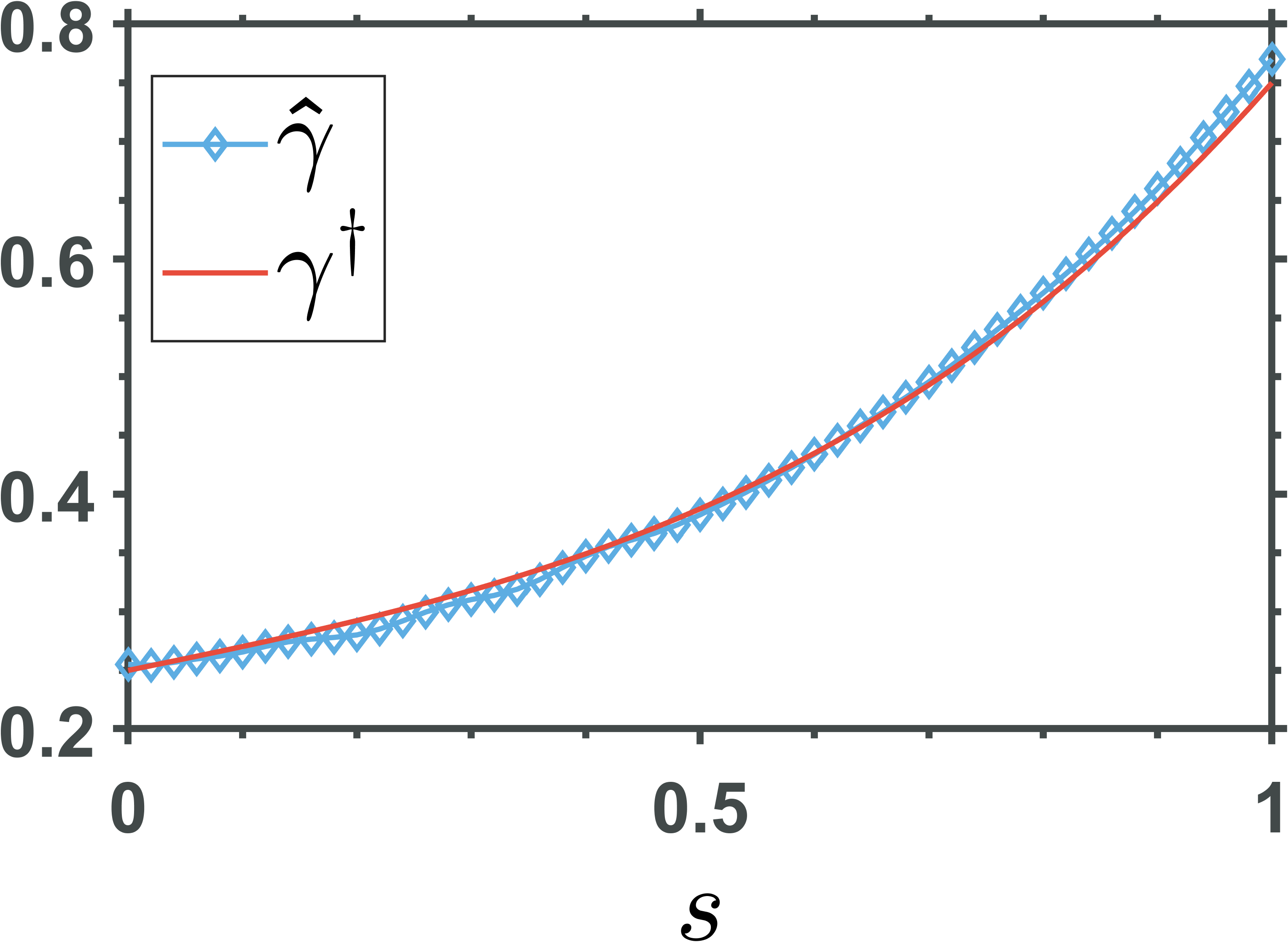} &
        \includegraphics[width=0.30\linewidth]{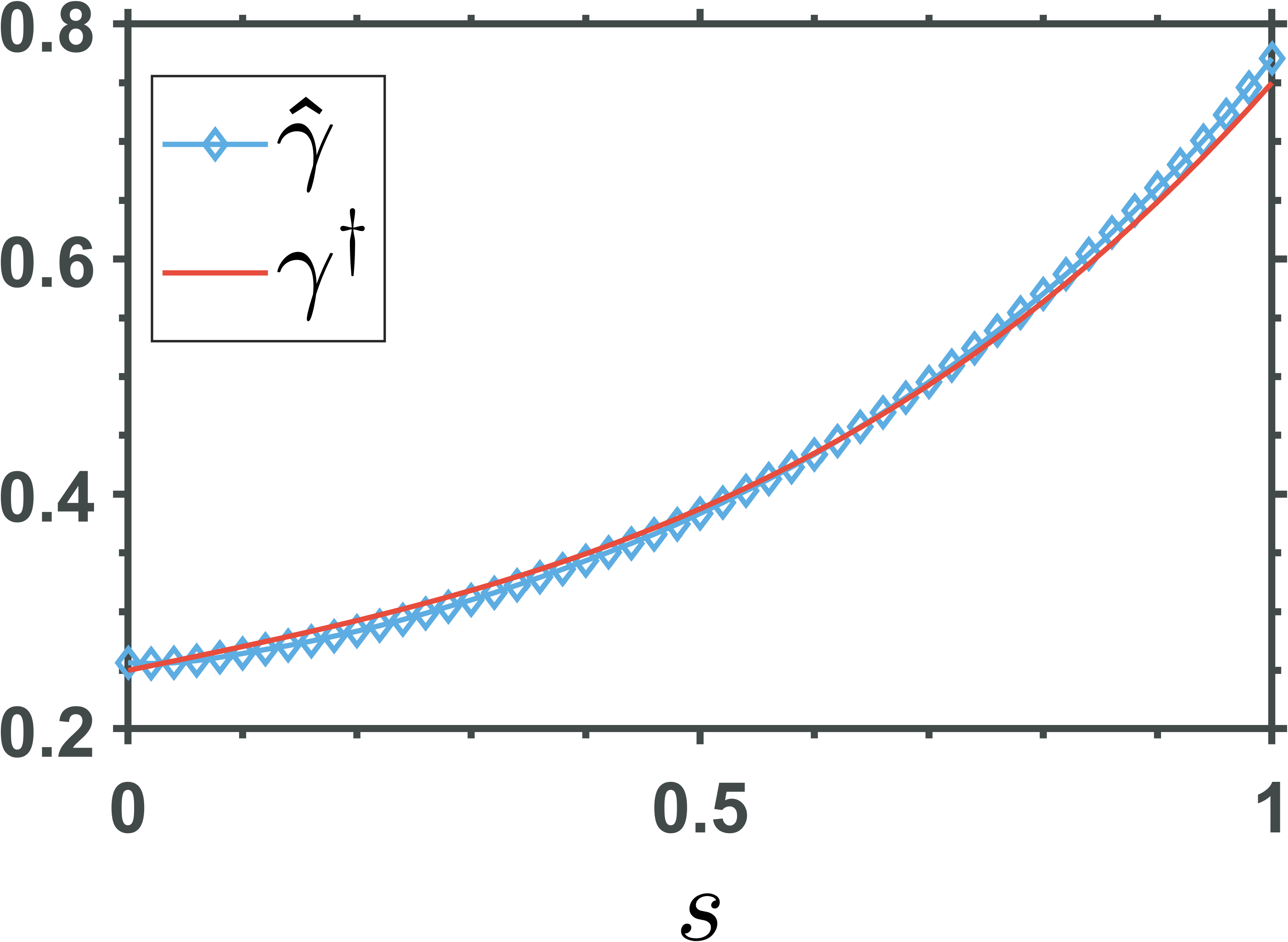} &
        \includegraphics[width=0.30\linewidth]{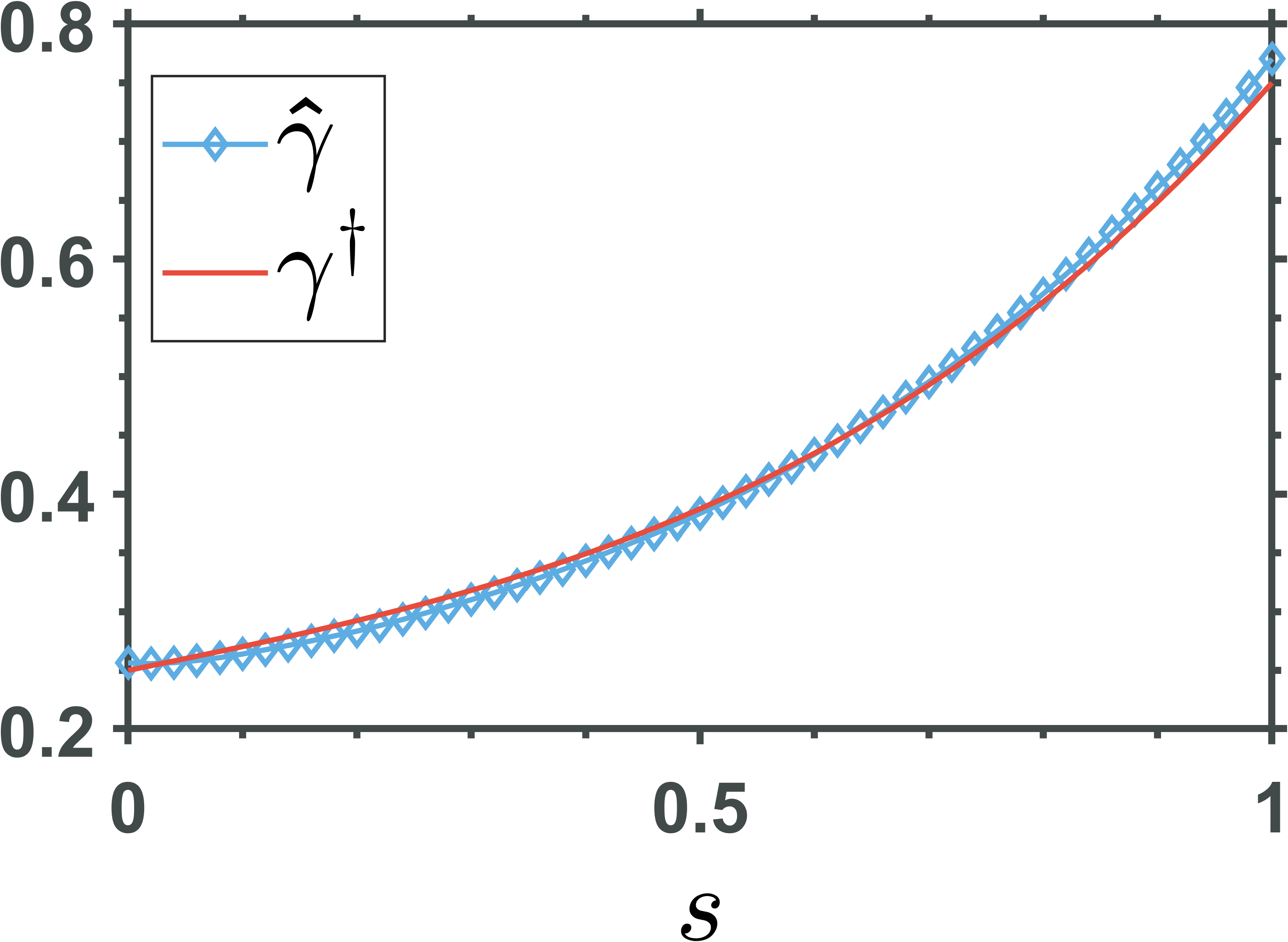} \\
          \includegraphics[width=0.30\linewidth]{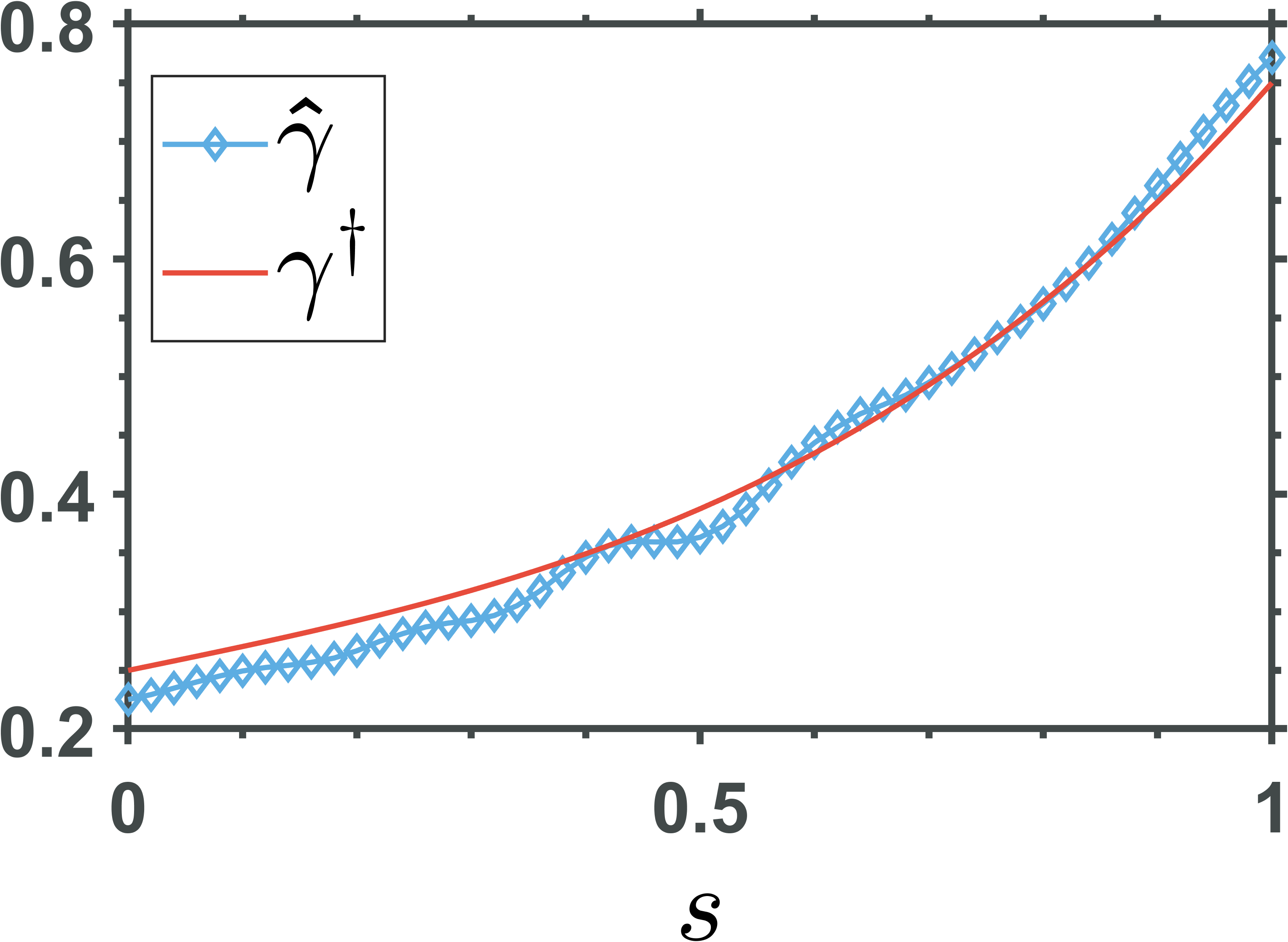} &
        \includegraphics[width=0.30\linewidth]{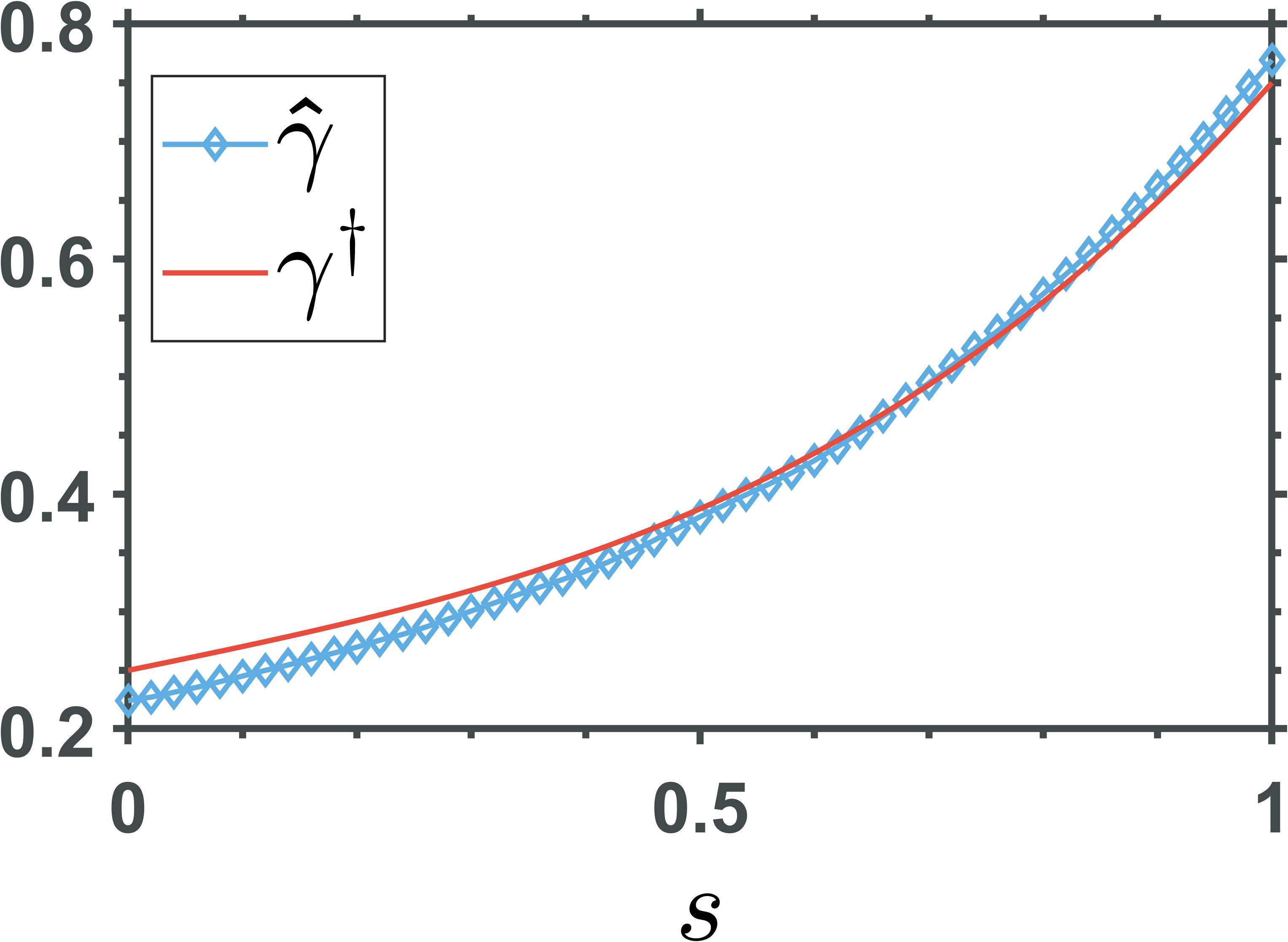} &
        \includegraphics[width=0.30\linewidth]{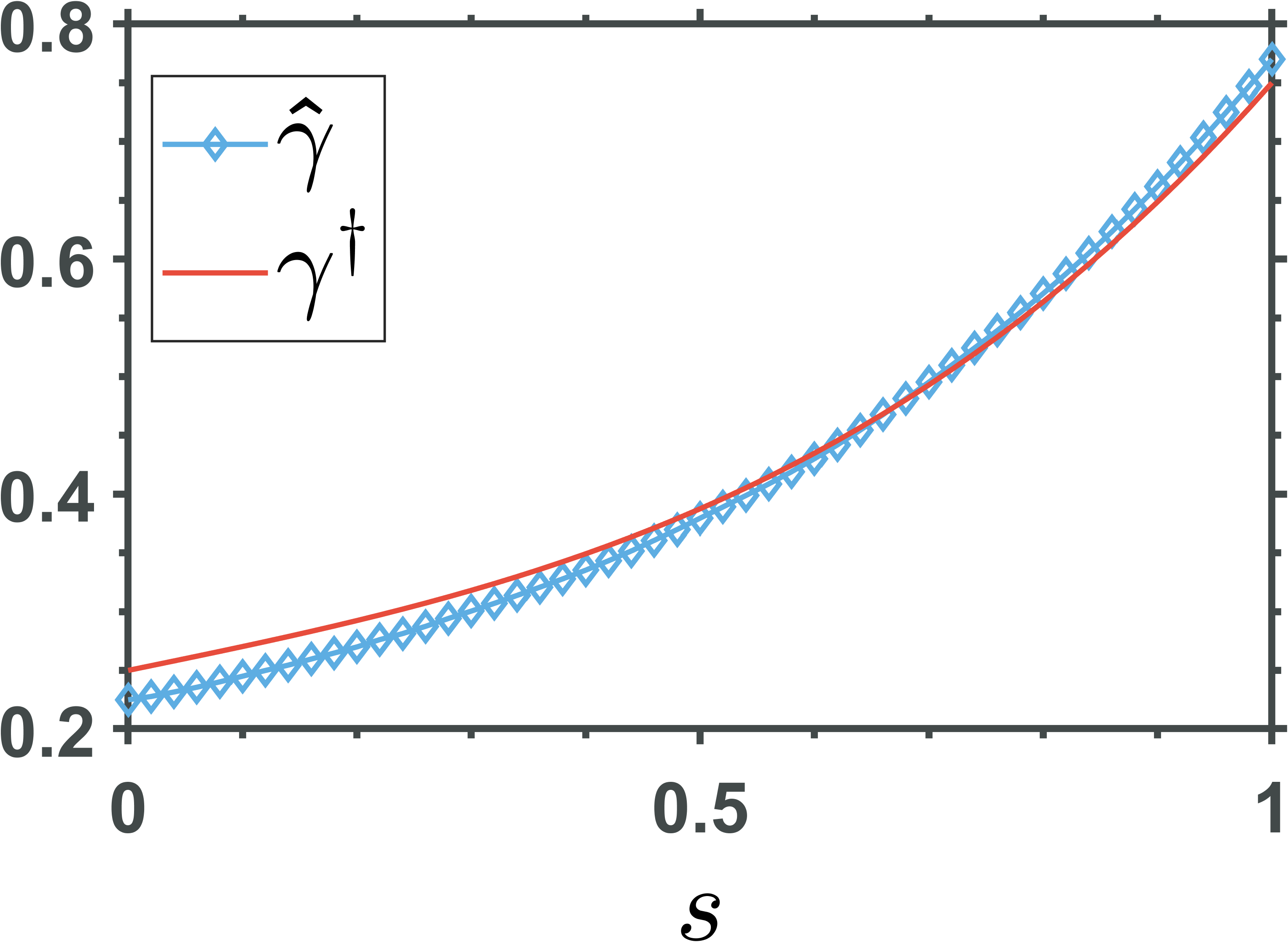} \\
        \includegraphics[width=0.30\linewidth]{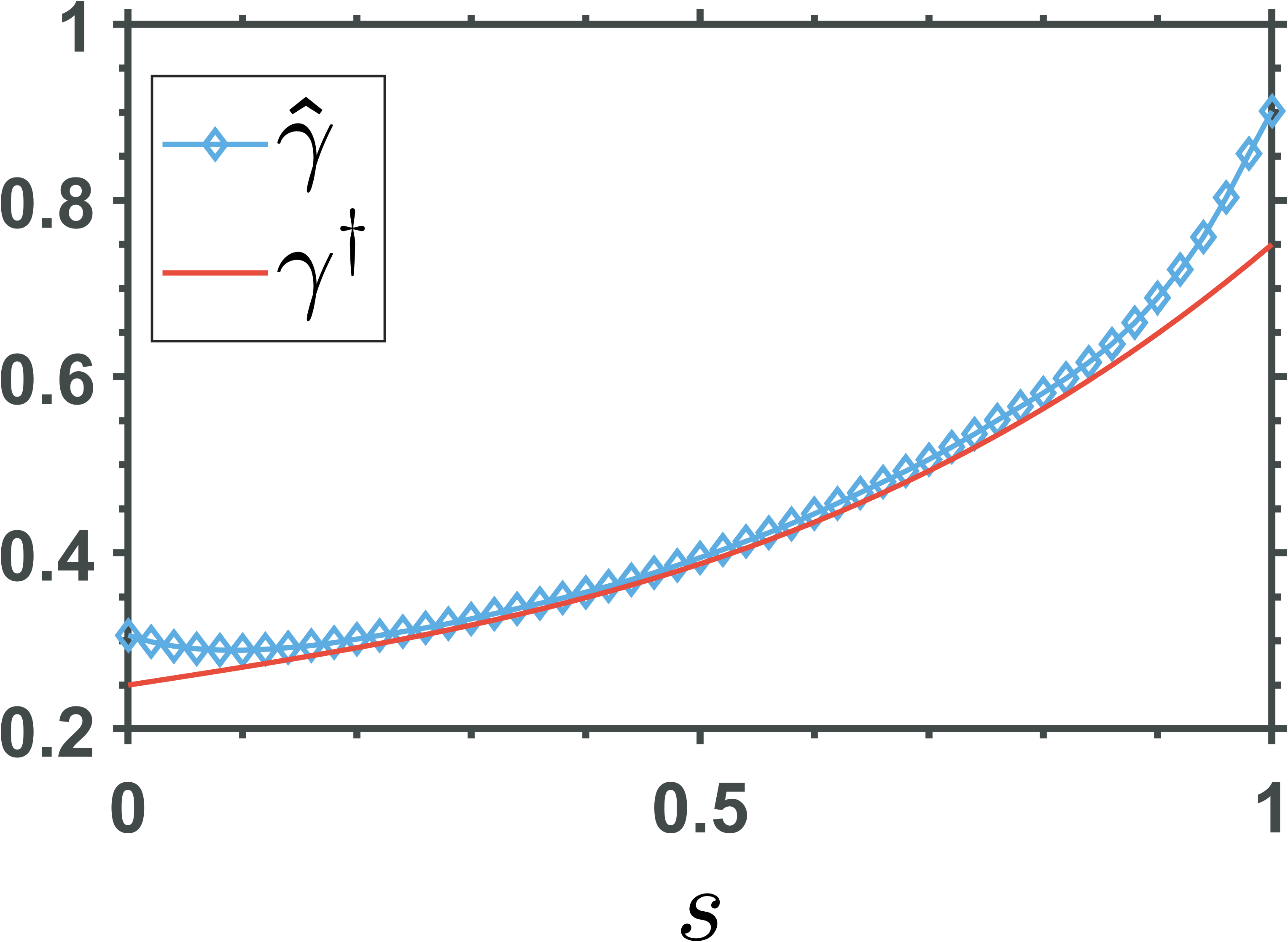} &
        \includegraphics[width=0.30\linewidth]{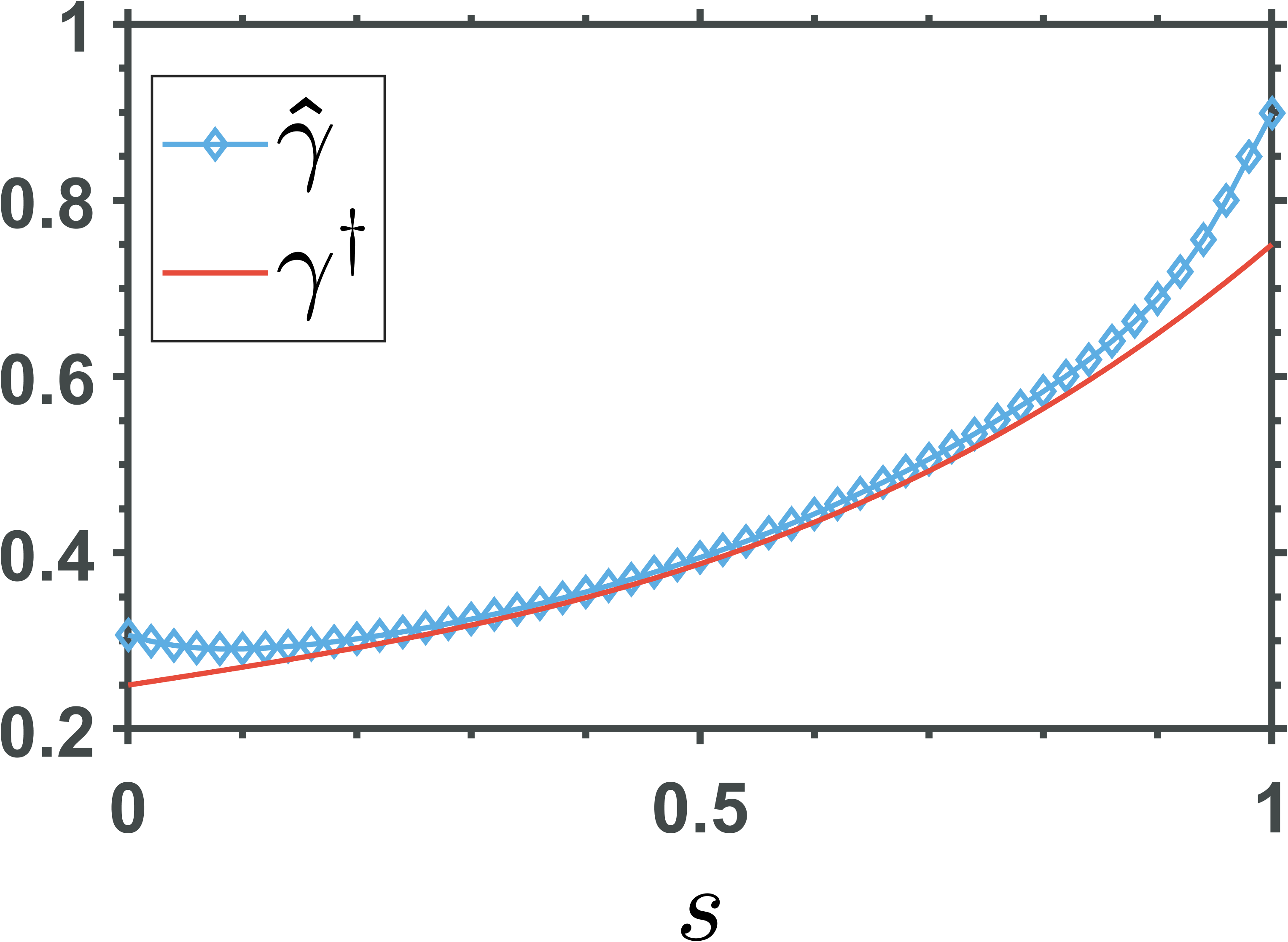} &
        \includegraphics[width=0.30\linewidth]{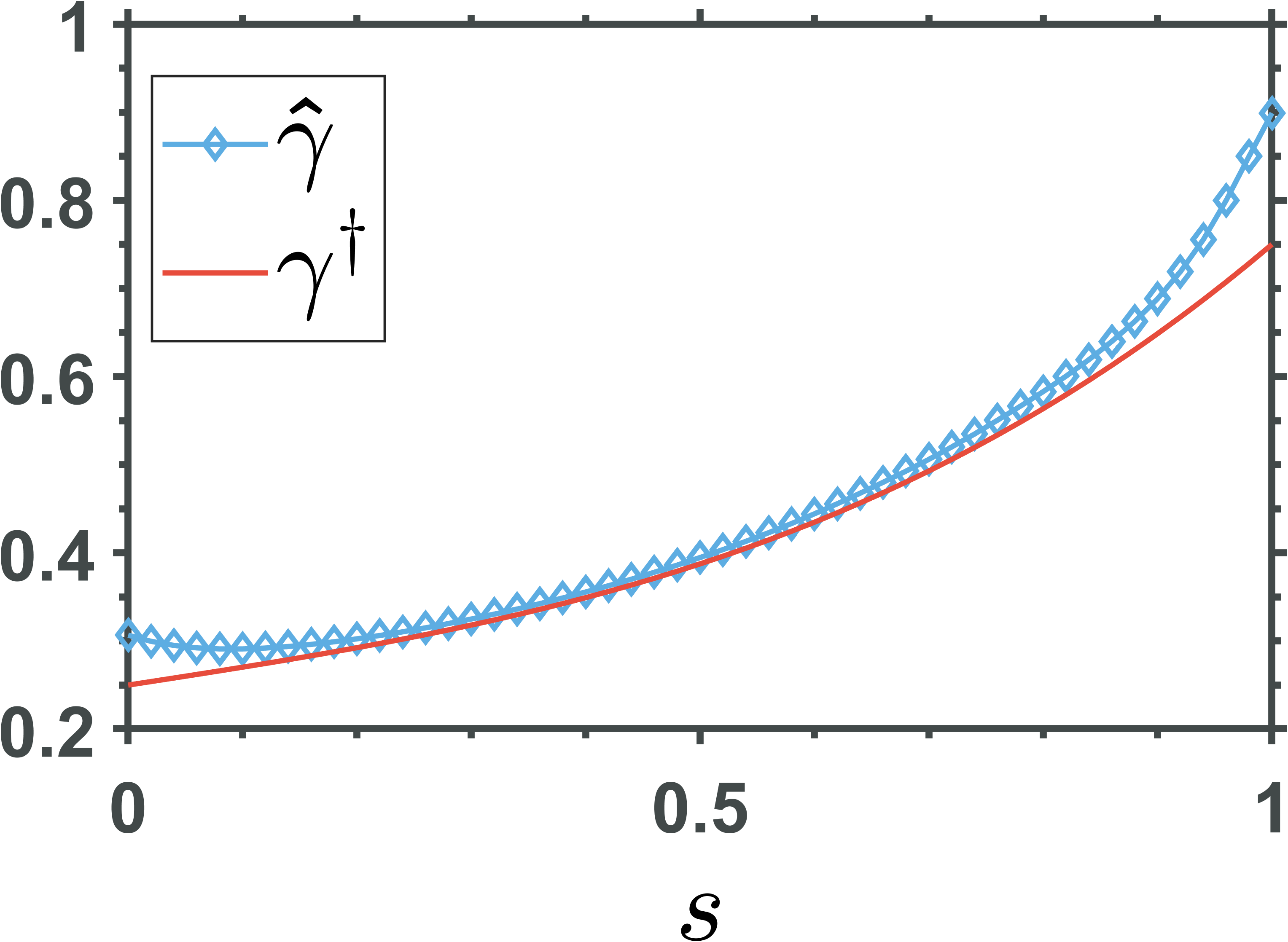} \\
         $\epsilon = 10^{-2}$ &  $\epsilon = 10^{-3}$ &  $\epsilon = 0$
    \end{tabular}
    \caption{
       The recovered conductivity $\hat{\gamma}$  versus the exact one $\gamma^\dag$ at various noise levels $\epsilon$. From top to bottom, the results are for $\Xi=\{1.1,1.2,\cdots,2.0\}$, $\Xi=\{2\}$ and $\Xi=\{1.1\}$, respectively. \\
    }
    \label{fig:gamma_reconstruction}
    \end{subfigure}
    \begin{subfigure}[b]{\linewidth}
        \centering
        \begin{tabular}{ccc}
       \includegraphics[width=0.30\linewidth]{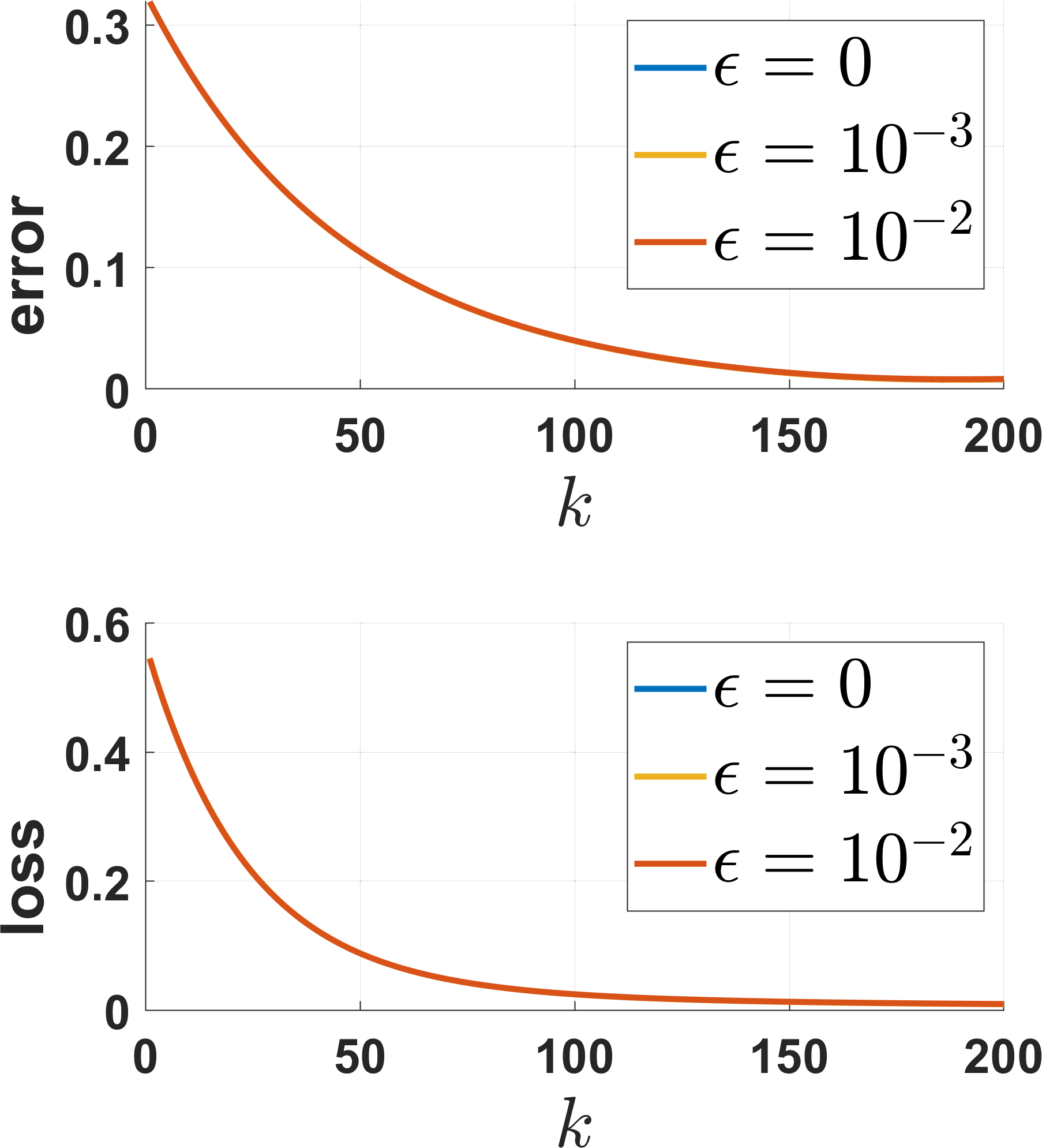}   & \includegraphics[width=0.30\linewidth]{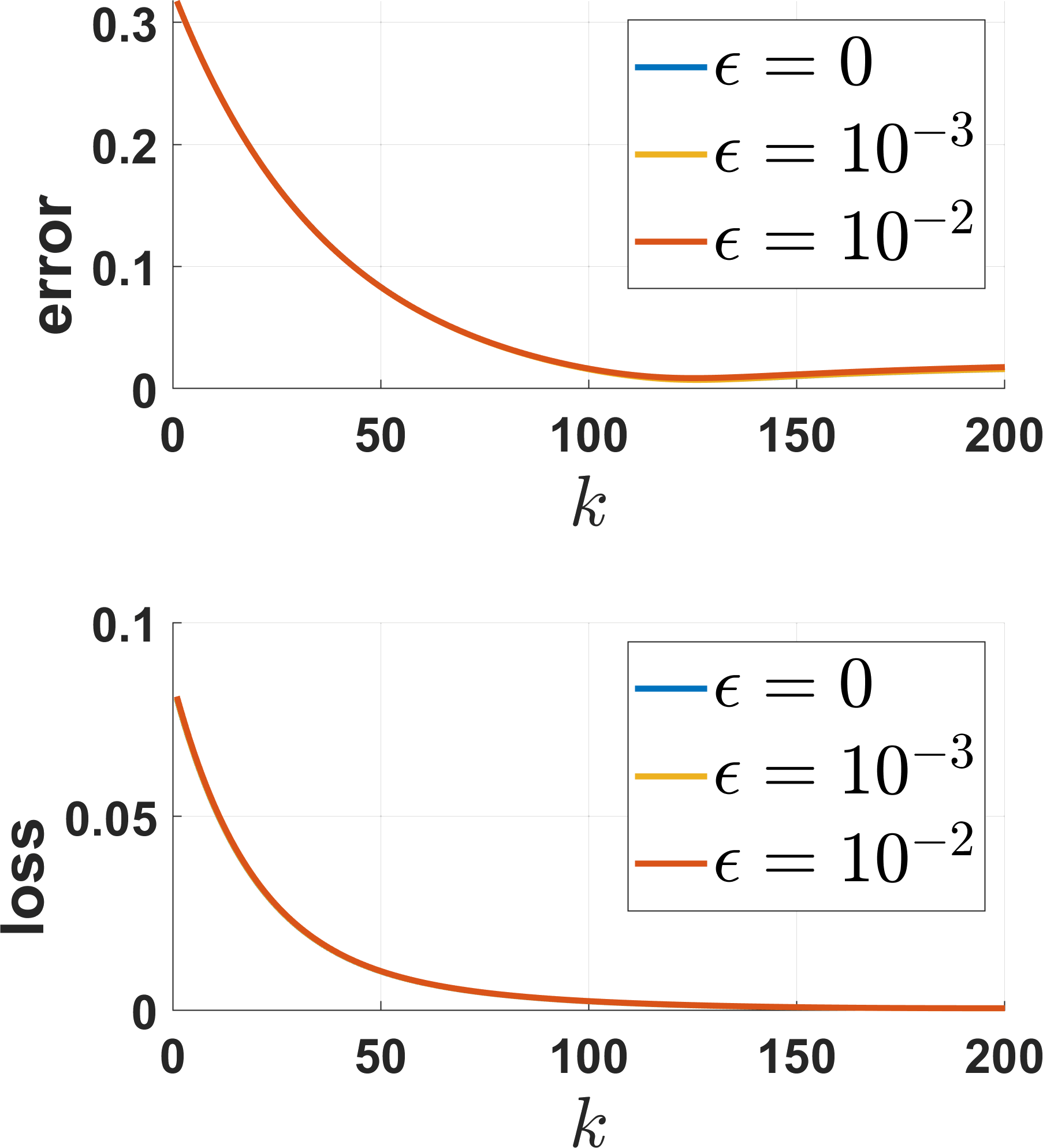}  &     \includegraphics[width=0.30\linewidth]{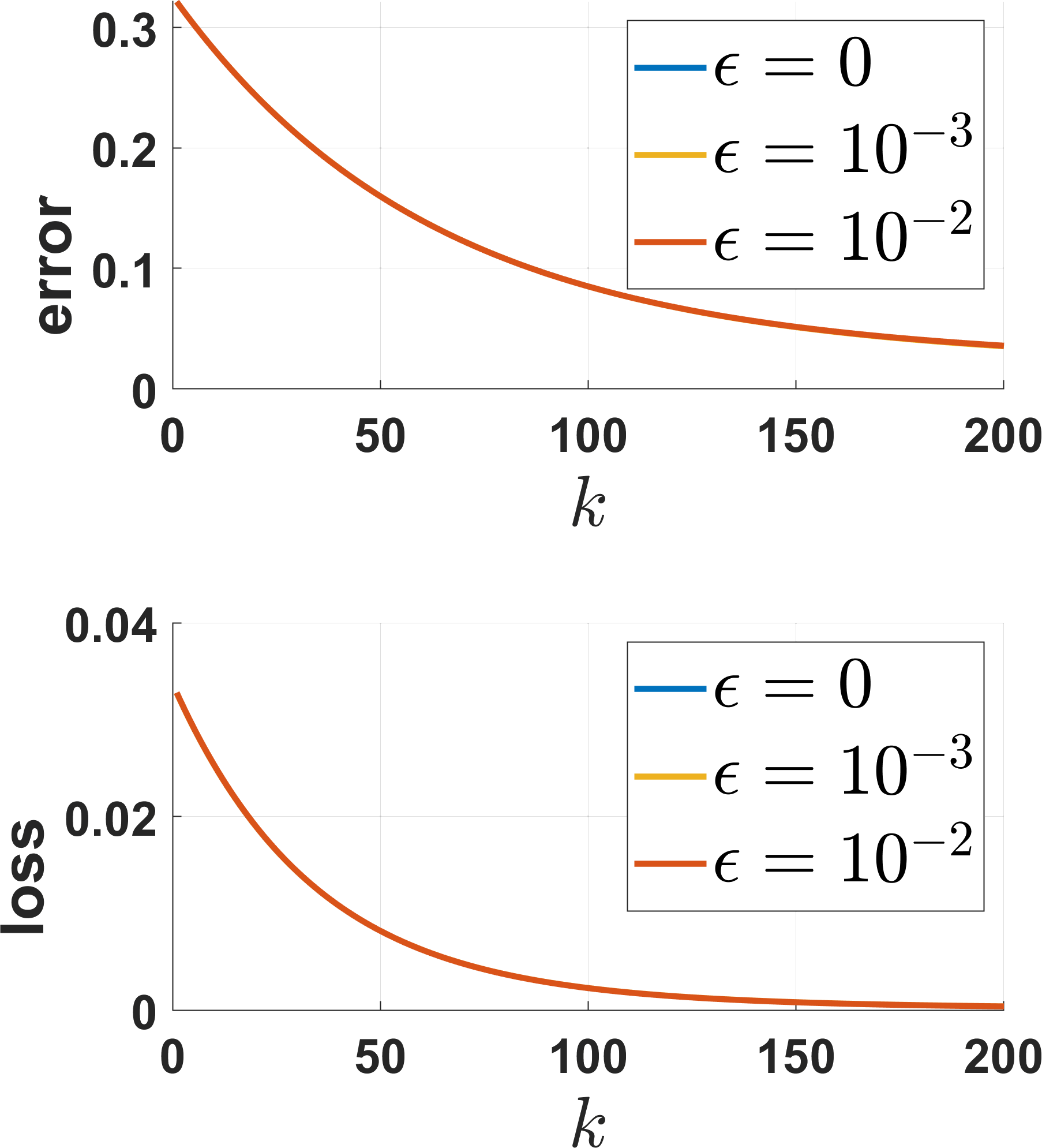}
\\
         $\Xi=\{1.1,1.2,\cdots,2.0\}$ &  $\Xi=\{2\}$ &  $\Xi=\{1.1\}$
    \end{tabular}
    \caption{The variation of the error $\|\hat\gamma-\gamma^\dag\|_{L^2((0,1))}$ (top) and the loss $ { {J}_0(\hat \gamma)}$ (bottom) with the iteration number $k$, for three different sets of measurement $\Xi$.  }
    \label{fig:Loss_Record}
    \end{subfigure}
    \caption{Numerical results for the 2D case.}
\end{figure}

First consider the 2D case on the unit disk \( \Omega = \mathbb{D} \) with the Dirichlet boundary condition(s) $g$ given by
$$
g_k(\theta) = k( e^{-\theta^2} - 0.1) , \quad (k,\theta) \in \Xi \times [-\pi,\pi],
$$
with \( \Xi = \{1 + 0.1d\}_{d=1}^{10} \). While \( g_k \) for \( d \geq 2 \) satisfies condition \eqref{eqn: cond on boundary fct. 1} in Theorem \ref{thm:main theoretical result}, the case \( d=1 \) violates the condition in the region near $s = 0.99$. We take the exact conductivity $\gamma^\dag$ to be \( \gamma^\dag(s) = 0.3s^2 + 0.2s + 0.25 \), and fix the regularization parameter $\beta$ at $10^{-3}$. The step size in the gradient descent method is set to $\frac{1}{|\Xi|}$, the initial guess is \( \gamma(s) = 0.5s + 0.5 \), sampled across 101 equidistant points in \( [-0.2,1.8] \). Fig. \ref{fig:gamma_reconstruction} presents the reconstruction results for both multi-measurement (\( \Xi = \{1.1:0.1:2.0\} \)) and single-measurement scenarios (\( \Xi = \{1.1\}, \{2.0\} \)) across three noise levels: \( \epsilon = 10^{-1} \), \( 10^{-2} \), and \( 10^{-3} \). For all three noise levels, the reconstruction $\hat\gamma$ for \( \Xi = \{2.0\} \) outperforms that for \( \Xi = \{1.1\} \) near \( s=1 \), due to the violation of crucial condition \eqref{eqn: cond on boundary fct. 1}. Moreover, the use of multiple measurements can significantly improve the stability of the reconstruction. The variation of the loss ${J}_0(\gamma)$ and the error $\|\hat\gamma-\gamma^\dag\|_{L^2((0,1))}$ during the optimization for the cases with three different noise levels are presented in Fig. \ref{fig:Loss_Record}. The difference of the training dynamics at different noise levels is very small, showing that the recovery is robust even in the presence of data noise.

\begin{figure}[hbt!]
    \centering
    \begin{subfigure}[b]{\linewidth}
        \centering
        \begin{tabular}{ccc}
            \includegraphics[scale=0.4]{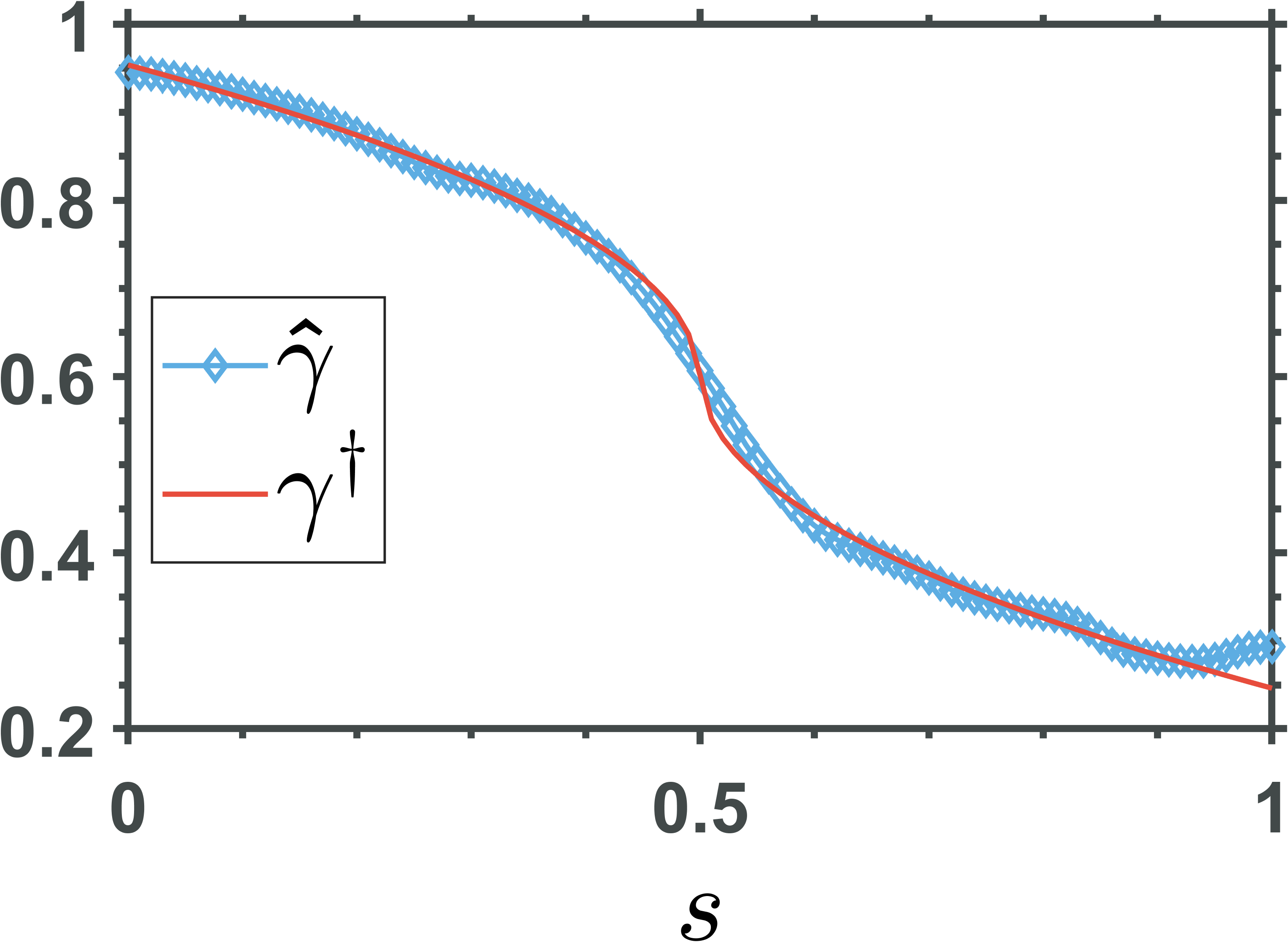} &
            \includegraphics[scale=0.4]{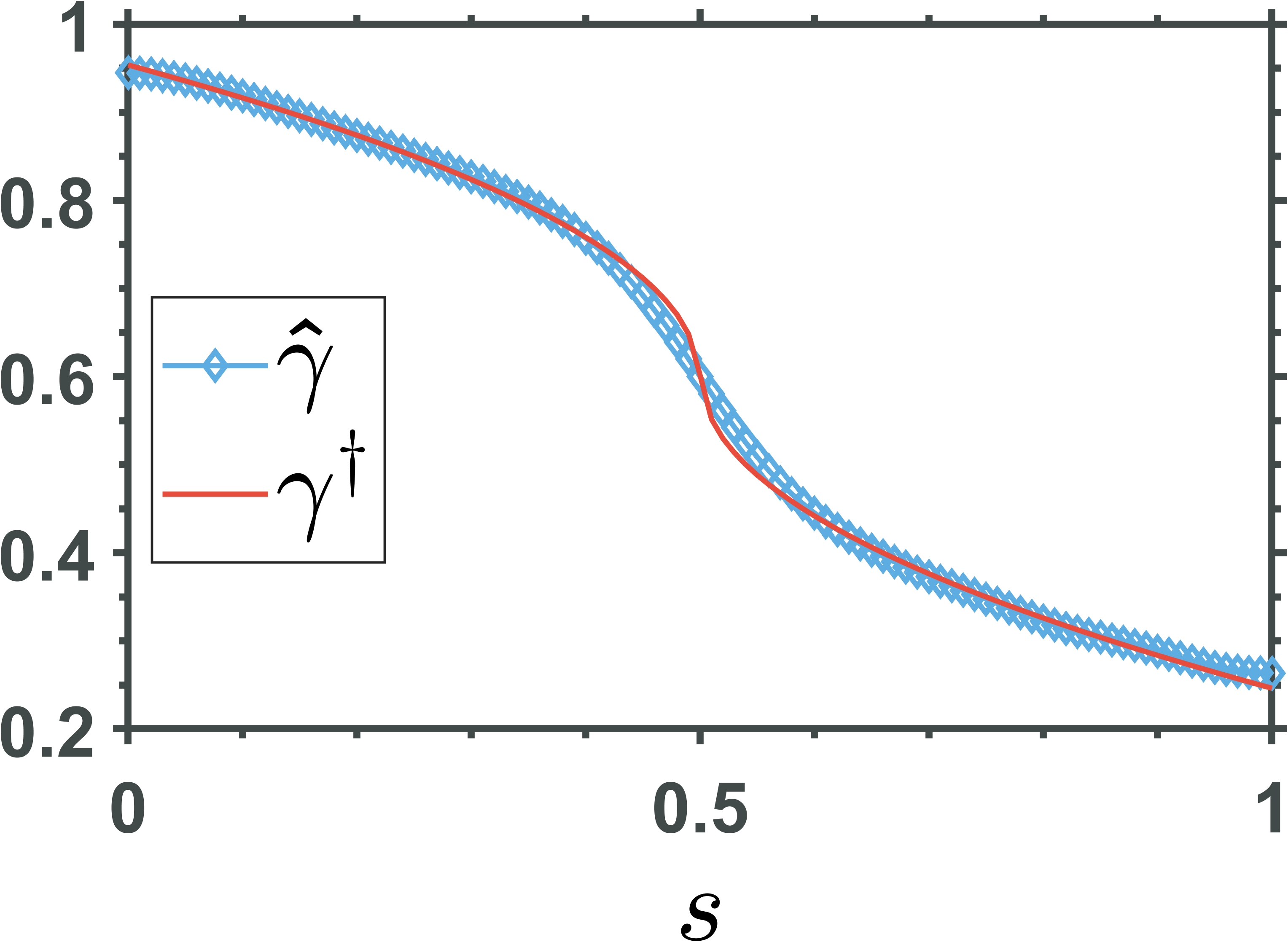} &
            \includegraphics[scale=0.4]{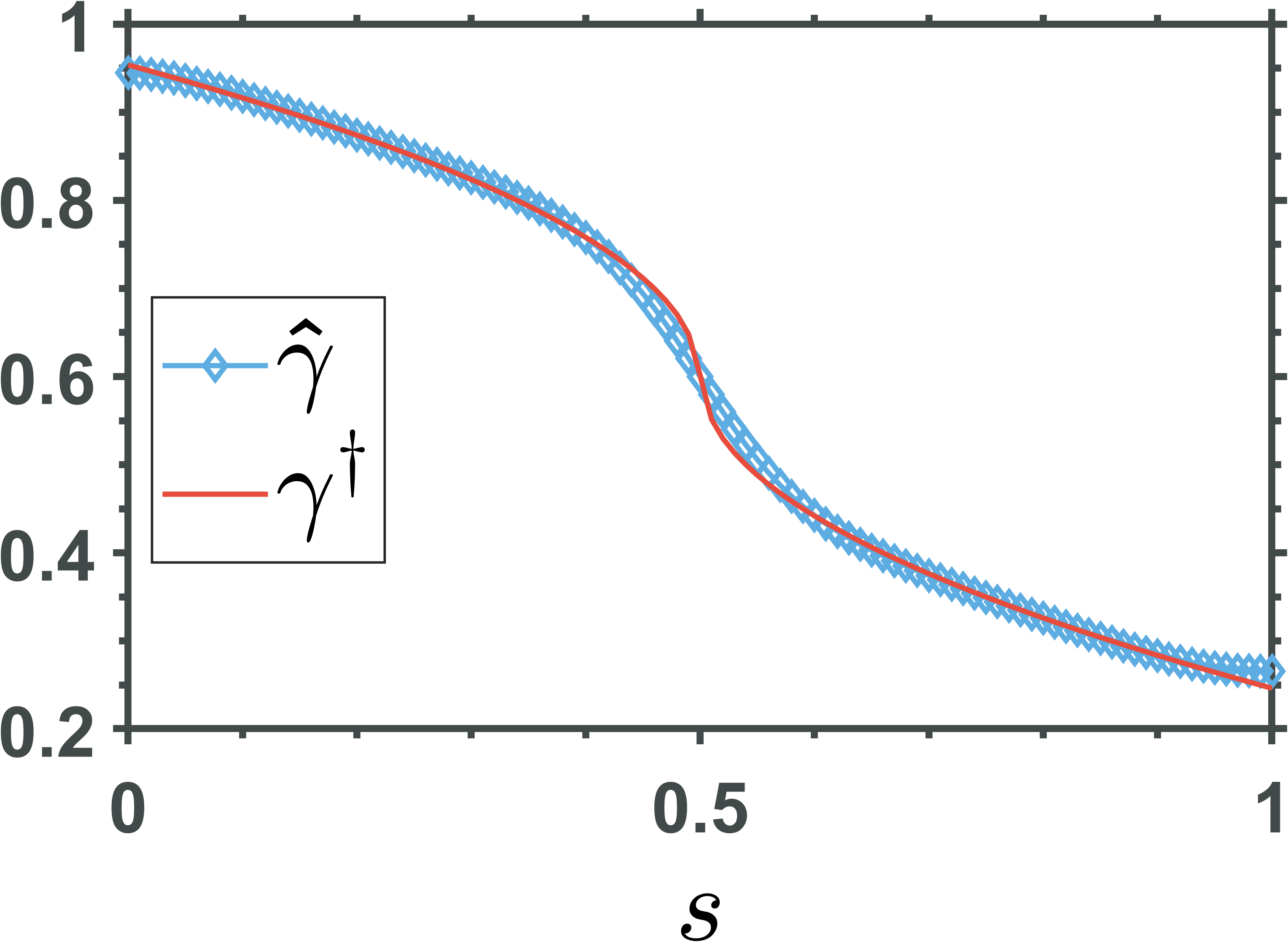} \\
        \end{tabular}
    \end{subfigure}
    \\
    \begin{subfigure}[b]{\linewidth}
        \centering
        \begin{tabular}{ccc}
            \includegraphics[scale=0.4]{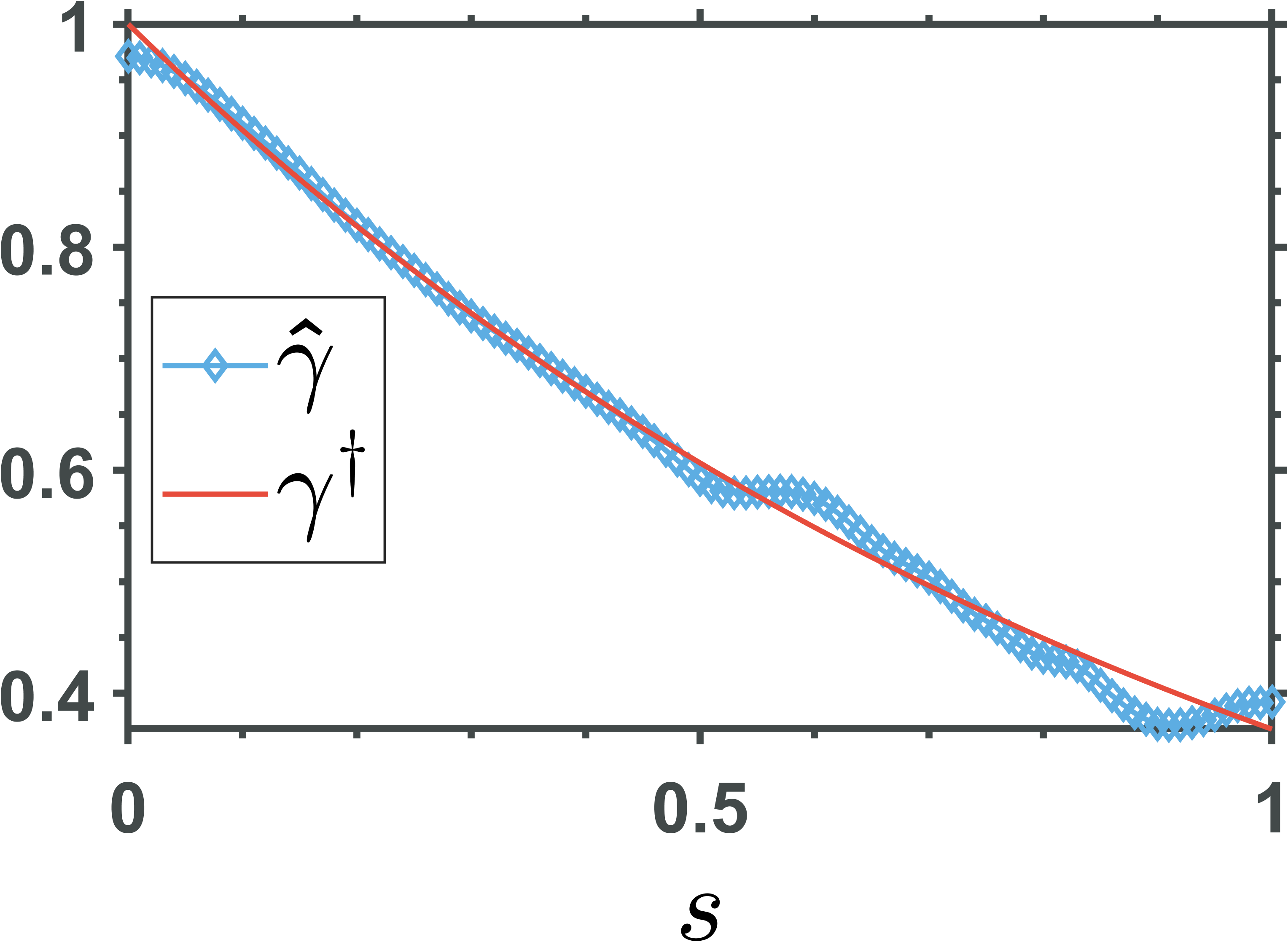}&
            \includegraphics[scale=0.4]{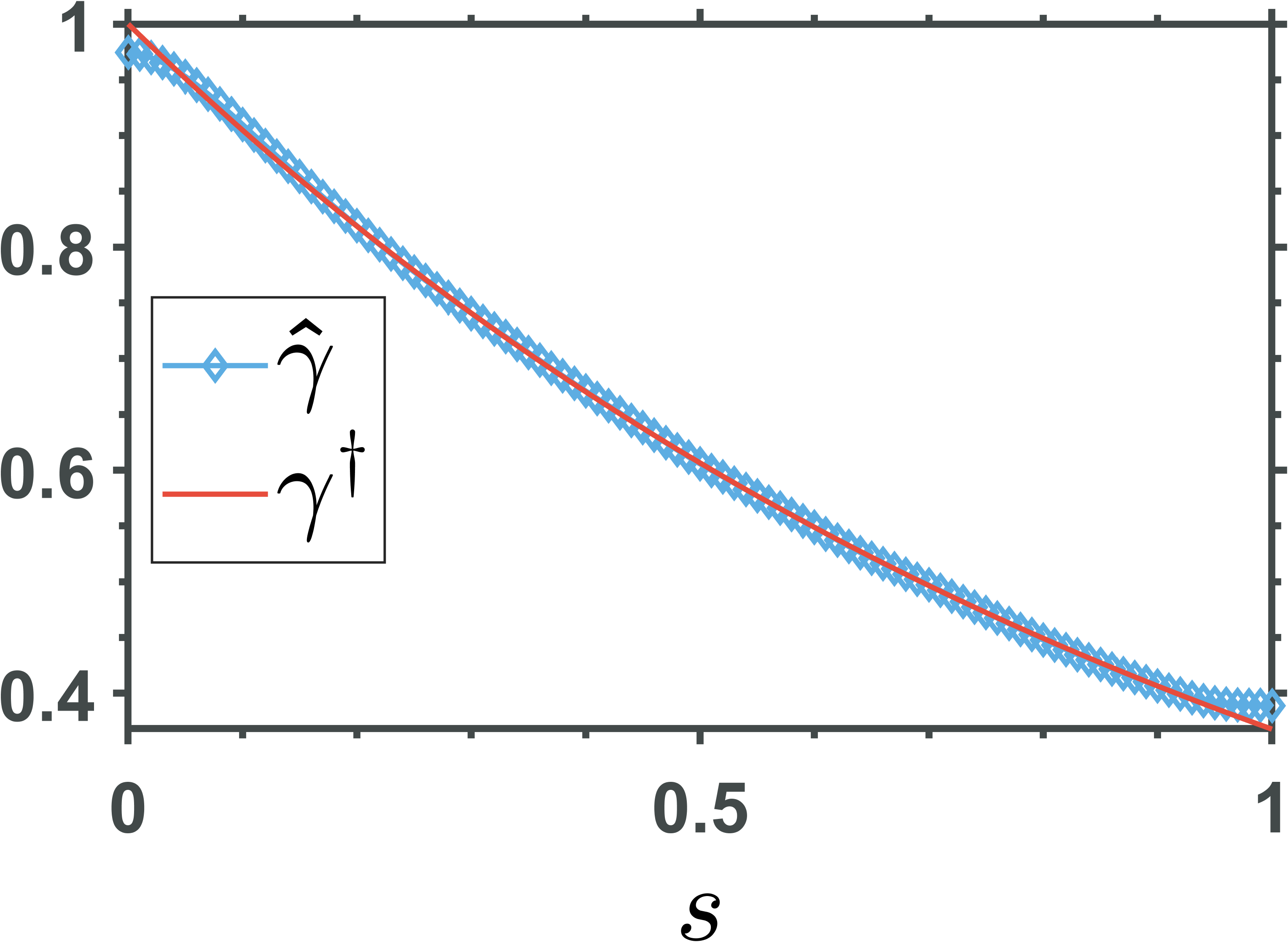} &
            \includegraphics[scale=0.4]{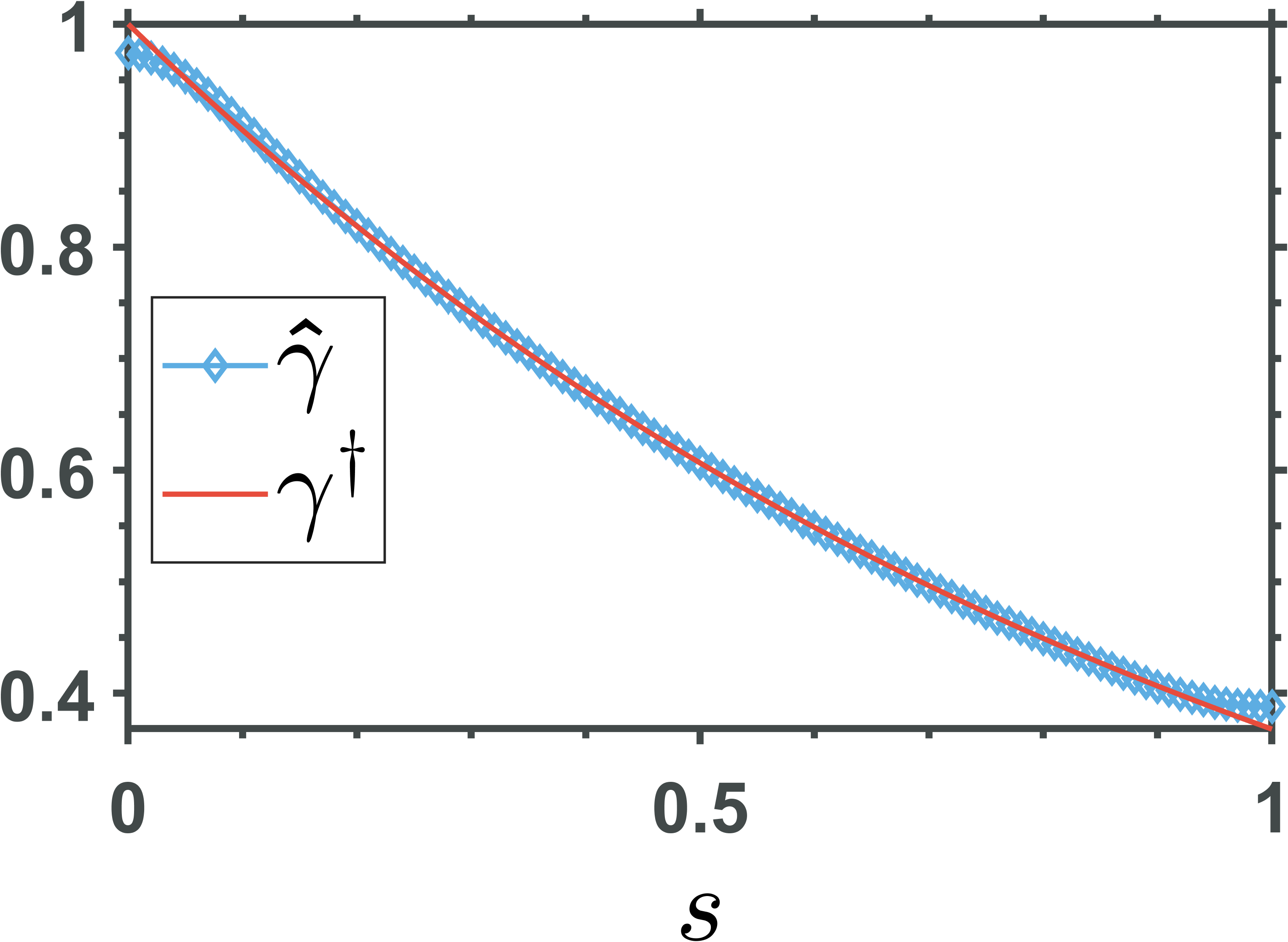} \\
            $\epsilon = 10^{-2}$ & $\epsilon = 10^{-3}$ & $\epsilon = 0$
        \end{tabular}
        \caption{ The recovered conductivities $\hat{\gamma}_1$ (top) and
        $\hat{\gamma}_2$ (bottom) at various noise levels $\epsilon$ versus the exact ones $\gamma_1^\dagger$ and $\gamma_2^\dagger$. }
    \end{subfigure}
    \\
    \begin{subfigure}[b]{\linewidth}
        \centering
        \begin{tabular}{cc}
        \includegraphics[width=0.485\linewidth]{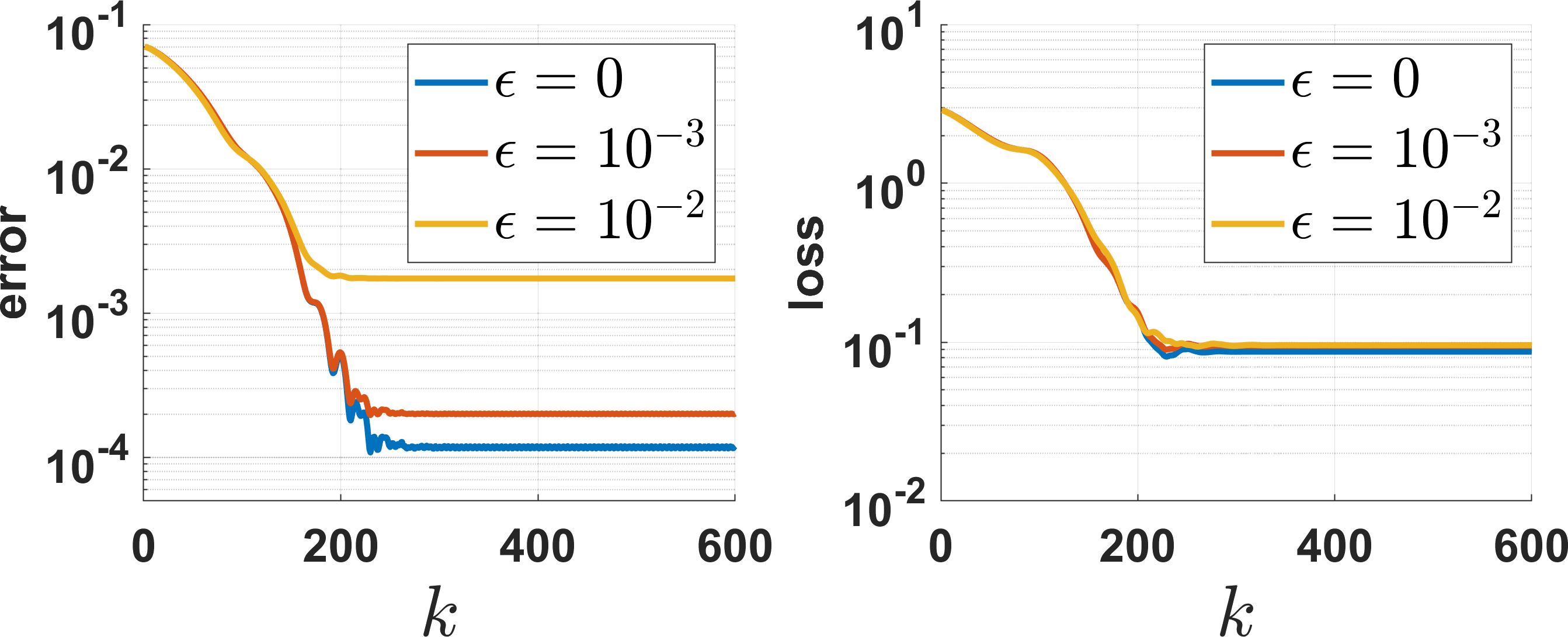} &
        \includegraphics[width=0.485\linewidth]{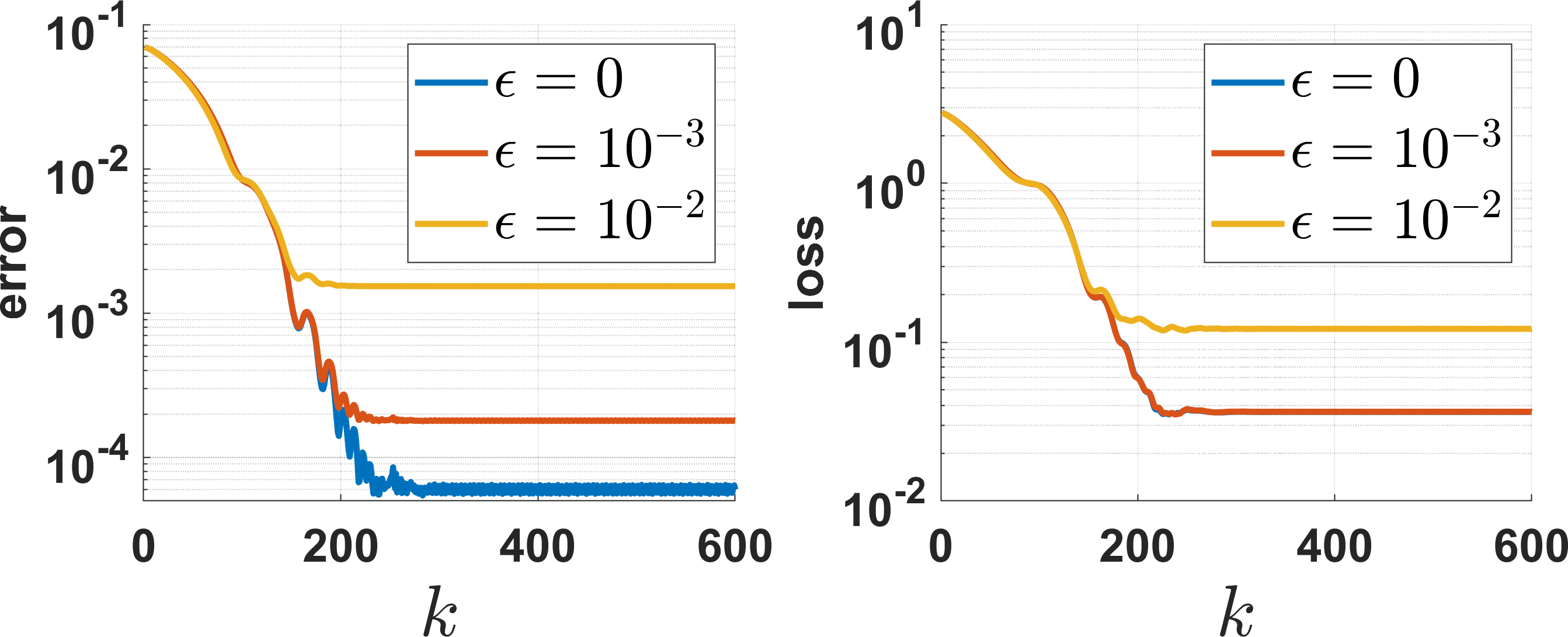}
        \end{tabular}
         \caption{The variation of the error  and the loss during the optimization for $\hat{\gamma}_1$ (left) and $\hat{\gamma}_2$ (right).  }
    \end{subfigure}
    \caption{Numerical results in the 1D case.}
    \label{fig:Loss_Record_dim1}
\end{figure}

Next consider the 1D case on the unit interval $\Omega=(0,1)$. We set boundary conditions \( u(0) = 0 \) and \( u(1) = \lambda \), with \( \lambda \) varying in the set \( \{ 0.01k : k = 1, \ldots, 100\} \), and examine two different cases for the conductivity \( \gamma \): a nonsmooth function \( \gamma_1^\dagger(s) = 1 - \operatorname{sgn}(s-0.5) \sqrt{|s-0.5|}\) and a smooth function \( \gamma_2^\dagger(s) = e^{-s} \). The regularization parameter $\beta$ is fixed at $\beta=0.1$. Due to the large number of nonlinear elliptic equations involved (100), we employ the Adam algorithm \cite{Kingma2014AdamAM} to minimize the functional \({J}_\beta \), with an initial guess of \( \gamma(s) = -0.25s + 0.5 \), step size $\frac{1}{300}$, full batch size 100 and both inertia momentum parameters set to $0.9$. The numerical results in Fig.~\ref{fig:Loss_Record_dim1} show comparable reconstruction accuracy for both conductivities. Thus, the regularity of \( \gamma^\dag \) does not affect much the reconstruction accuracy. Moreover, the reconstruction at the noise level $\varepsilon = 10^{-3}$ is comparable with that for the exact data $v$, indicating the excellent  stability of the reconstruction. The reconstruction remains accurate near the boundary points $s=0$ and $s=1$, which agrees with the prediction by Theorem \ref{thm:1D-stability}.

\section{Conclusion}
In this work we have established several novel H\"{o}lder type stability results for the inverse problem of recovering a quasilinear term in an elliptic equation from the conormal data on the
 boundary in both one- and multi-dimensional cases, and have presented several numerical tests illustrating the feasibility of the stable reconstruction. It is of much interest to utilize the stability results for the numerical analysis of the regularized reconstructions and the discrete approximations (e.g., with the finite element method). Theoretically, it is of much interest to analyze the optimality of the H\"{o}lder stability in the multi-dimensional case.

\appendix
\section{Technical results for the proof of Theorem \ref{thm:1D-stability}}
\label{sec:appendix}
\begin{figure}[hbt!]
        \centering
        \includegraphics[width=0.5\linewidth]{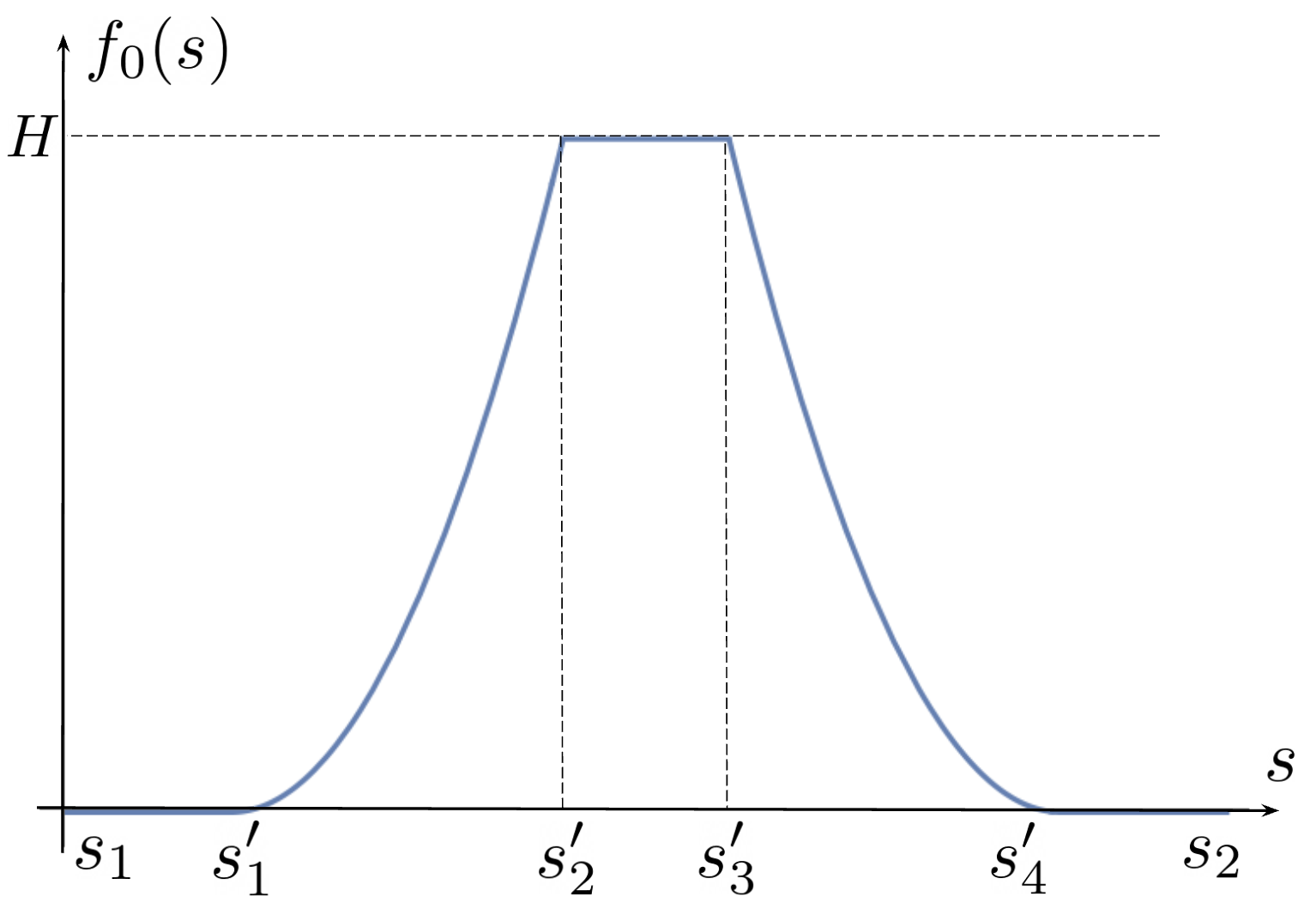}
       \caption{A schematic illustration of the minimizer $f_0(s)$.}
        \label{fig:Possible form 1}
\end{figure}

In the appendix we collect several technical results that are used in the proof of Theorem \ref{thm:1D-stability}. The first result gives an explicit form of the minimizer $f_0$ under the additional restriction that $f(s_1)=f(s_2)=0$.
\begin{proposition}
\label{prop:Needyform1}
Suppose $I=(s_1,s_2)$ and denote
$$
\overline{\mathfrak{M}}_{I,H,S}^p = \{ f \in \mathfrak{M}_{I,H,S}^p \mid f(s_1) = f(s_2) = 0 \},
$$
with $\mathfrak{M}_{I,H,S}^p$ defined in \eqref{eqn: defn of M^p_H,S}. Let $f_0 = \arg \min_{f \in \overline{\mathfrak{M}}_{I,H,S}^p} \mathcal{F}^p_{I}(f)$, with $\mathcal{F}_{I}^p$ defined in \eqref{eqn:defn of F_I^p}. Then $f_0$ must take the following form {\rm(}see Fig. \ref{fig:Possible form 1} for a schematic illustration{\rm)}: for some  $s_1\leq s_1'\leq s_2' \leq s_3' \leq s'_4 \leq s_2$ and some constants $b_1,b_2,\tilde{b}_1,\tilde{b}_2,c_1$ and $c_2$,
    \begin{equation}
    \label{eqn:Needyform1}
    f_0(s) = \begin{cases}
     c_1\left(s+b_1\right)^{\frac{p}{p-1}}+\tilde{b}_1 , & s_1' \leq s < s_2', \\
       H, & s_2' \leq s < s_3', \\
      c_2\left(s+b_2\right)^{\frac{p}{p-1}}+\tilde{b}_2 , & s_3' \leq s < s_4', \\
      0, & \text{otherwise}.
    \end{cases}
\end{equation}
\end{proposition}
\begin{proof}
Let $J = \{s\in(s_1,s_2):f_0(s) = H \}$, $\tilde{s}_1=\inf J$ and $\tilde{s}_2=\sup J$. We show that the set $J$ is connected. Suppose $f_0\not\equiv H$ on $[\tilde{s}_1,\tilde{s}_2]$. Let $\beta_1 = H^{-1}\int_{\tilde{s}_1}^{\tilde{s}_2}f_0(s) \d s$. Then we construct $\tilde{f}_0$ as
$$
    \tilde{f}_0(s) =
\begin{cases}
f_0(s), & s_1 \leq s < \tilde{s}_1 , \\
H, & \tilde{s}_1  \leq s < \beta_1+ \tilde{s}_1 , \\
f_0(s+\tilde{s}_2-\beta_1-\tilde{s}_1), & \beta_1+\tilde{s}_1  \leq s < s_2-\tilde{s}_2+\beta_1+\tilde{s}_1 , \\
0,& \text{otherwise}.
\end{cases}
$$
We can check that $\tilde{f}_0 \in \overline{\mathfrak{M}}_{I,H,S}^p$ and we also have
$$
\|f_0\|_{L^1(I)} = \|\tilde{f}_0\|_{L^1(I)}\quad\mbox{and} \quad \mathcal{F}^p_I(f_0) - \mathcal{F}^p_I(\tilde{f}_0) = \int_{\tilde{s}_1}^{\tilde{s}_2} [f_0'(s)]^p \d s > 0,
$$
which contradicts the assumption that $f_0$ is a minimizer. This proves that $f_0 \equiv H$ on $[\tilde{s}_1,\tilde{s}_2]$. Next, we prove the connectivity of the set $K=\{s:f_0(s) > 0\}$. Suppose that $K =  K_0 \cup K'$, with $J \subset K_0$, $K' \neq \emptyset$ and $K_0 \cap K' = \emptyset$.
Let $K_0 = [\hat{s}_1,\hat{s}_2]$ and $\beta_2 = H^{-1} \int_{K'} f_0(s) \d s $. Then we construct $\hat{f}_0$ as
$$
\hat{f}_0(s) =
\begin{cases}
    f_0(s+\hat{s}_1-s_1 ), & s_1 \leq s <\tilde{s}_1-\hat{s}_1+s_1 , \\
    H, & \tilde{s}_1-\hat{s}_1 +s_1 \leq s < \tilde{s}_2 + \beta_2-\hat{s}_1 +s_1, \\
    f_0(s-\beta_2+\hat{s}_1-s_1 ), & \tilde{s}_2+\beta_2-\hat{s}_1+s_1 \leq s < \hat{s}_2 + \beta_2-\hat{s}_1+s_1, \\
    0, & \text{otherwise}. \\
\end{cases}
$$
We similarly have $\hat{f}_0 \in \overline{\mathfrak{M}}_{I,H,S}^p$ and
$$
\|f_0 \|_{L^1(I)} = \|\hat{f}_0\|_{L^1(I)}\quad \mbox{and} \quad \mathcal{F}^p_I(f_0) - \mathcal{F}_I^p(\hat{f}_0) = \int_{K'} [f'_0(s)]^p \d s >0 ,
$$
which again contradicts the assumption that $f_0$ is a minimizer. Hence the set $K$ is connected. By repeating the argument in the proof of Theorem \ref{thm:1D-stability}, on $K\setminus J$, we deduce that $f_0$ must take the form
$$
f_0(s)= c(s+b)^{\frac{p}{p-1}} + \tilde{b}.
$$
By combining this with the fact that both $J$ and $K$ are connected, we directly deduce that $f_0$ must be of the form \eqref{eqn:Needyform1}, which concludes the proof of the proposition.
\end{proof}

The second result similarly gives an explicit form for $f_0$ under the condition that $f(s_1) = 0$ and $s_2 = R$. The proof is analogous to that of Proposition \ref{prop:Needyform1} and is thus omitted.
\begin{proposition}
  \label{prop:Needyform2}
   Suppose $I=(s_1,R)$ and denote
    $$
    \widetilde{\mathfrak{M}}_{I,H,S}^p = \{ f \in \mathfrak{M}_{I,H,S}^p \mid f(s_1) = 0 \},
    $$
    with $\mathfrak{M}_{I,H,S}^p$  defined in \eqref{eqn: defn of M^p_H,S}. Let $f_0 = \arg \min_{f \in \widetilde{\mathfrak{M}}_{I,H,S}^p} \mathcal{F}^p_{I}(f)$, with  $\mathcal{F}_{I}^p$ defined by \eqref{eqn:defn of F_I^p}. Then $f_0$ must take the following form: for some  $s_1 \leq s'_1 \leq s'_2 \leq R$ and some constants $b$, $\tilde{b}$ and $c$,
\begin{equation}
   f_0(s) = \begin{cases}

        c(s+b)^{\frac{p}{p-1}}+\tilde{b}, & s'_1 \leq s < s'_2, \\
        H, & s'_2 \leq s < R, \\
        0, & \text{otherwise}.
    \end{cases}\label{eqn:Needyform2}
\end{equation}
\end{proposition}

The last result gives an estimate of the minimal value $\mathcal{F}^p_I(f_0)$.

\begin{proposition}
\label{prop:minimal_F}
Let $f_0\in \mathfrak{M}_{I,H,S}^p$ be defined by either \eqref{eqn:Needyform1} or \eqref{eqn:Needyform2}. Then, we have
    \begin{equation}
    \label{eqn:estimate for F_I^p}
\mathcal{F}_I^p(f_0) \geq \left( \frac{p-1}{2p-1} \right)^{p-1} \frac{H^{2p-1}}{S^{p-1}}.
 \end{equation}
\end{proposition}
\begin{proof}
The proof is lengthy and we split it into three steps.
In the first two steps we prove the estimate \eqref{eqn:estimate for F_I^p} when $f_0$ is of the form \eqref{eqn:Needyform2}, i.e.,
\begin{equation}\label{eqn:f0-1} f_0(s)=c(s+b)^{\frac{p}{p-1}}+\tilde b,\quad \forall s\in (s_1',s_2')\subset [0,R],
\end{equation}
with $f_0(s_1')=0$ and $f_0(s_2')=H$, and $f_0(s)>0$ over $s\in (s_1',s_2')$.\\
\noindent\textbf{Step 1: Estimate $\mathcal{F}_I^p(f_0)$ by $\ell:=s_2'-s_1'$.} We first determine conditions on the constants $b$, $\tilde b$ and $c$. The condition $ f_0(s'_1) = 0$ implies $\tilde{b} = -c(s_1'+b)^{\frac{p}{p-1}}$. Since $\tilde{f}_0(s) > 0 $ for $s \in (s_1',s_2']$, we deduce
\begin{equation}\label{eqn:f0-comp}
c(s+b)^{\frac{p}{p-1}} > c(s_1'+b)^{\frac{p}{p-1}}, \quad \forall s \in (s_1',s_2'].
\end{equation}
Clearly $c\neq0$. Thus, for $s>s_1'$ close to $s_1'$, we have $f_0'(s_1')=\frac{cp}{p-1}(s_1'+b)^{\frac{1}{p-1}}\geq0$. Then there are two possible cases: (i) $b\geq -s_1'$ and $c>0$ and  and (ii) $c<0$ and $b< -s_1'$. However, in case (ii), for $s>s_1$ and close to $s_1$, the inequality \eqref{eqn:f0-comp} would fail to hold.
Thus we must have $c>0$, $b\geq -s_1'$ and $\tilde b = -c(s_1'+b)^{\frac{p}{p-1}}$ in \eqref{eqn:f0-1}, and we have $f_0'(s)= \frac{p}{p-1} c(s+b)^{\frac{1}{p-1}}\geq 0$ on $(s_1',s_2')$ and clearly $f_0'(s) = 0$ on $(s_1,s_1')\cup(s_2',R)$. Below let $\tilde f_0(s)= f_0(s+s_1')$ for $s\in (0,\ell)$, which is a translated version of $f_0$, and $d=b+s_1'$. Thus,
    $$
    \int_{s_1}^{R} |f_0'(s)| \d s = \int_{0}^{\ell}  \tilde{f}_0'(s) \d s = \tilde{f}_0(\ell)-\tilde{f}_0(0) = H.
    $$
Then by H\"older's inequality, we have
    $$
    H^p = \left(\int_{s_1}^{R} |f_0'(s)| \d s\right)^p \leq \ell^{p-1} \int_{0}^{\ell} |\tilde{f}_0'(s)|^p \d s  =\ell^{p-1}\mathcal{F}_I^p(f_0).
    $$
It follows directly that
\begin{equation}
\label{eqn:ineq in previous A.3}
    \mathcal{F}_I^p(f_0) \geq \frac{H^p}{\ell^{p-1}}.
\end{equation}
\textbf{Step 2: Estimating $\ell$ by $S/H$.} By the fundamental theorem of calculus, we have $
H = \int_{0}^{\ell} \tilde{f}_0'(s) \d s .$
Moreover, let
$\widetilde S = \int_{0}^{\ell} \tilde{f}_0(s) \d s.$
By changing the order of integration, we obtain
$$
  \int_{0}^{\ell} (\ell - s)\tilde{f}_0'(s) \d s = \int_{0}^{\ell} \int_{0}^s \tilde{f}_0'(\tau) \d\tau \d s = \int_{0}^{\ell} \tilde{f}_0(s) \d s = \widetilde S.
$$
By combining the identities for $H$ and $S_1$ and \eqref{prop:Needyform2}, we derive
\begin{equation}
\label{eqn: key idenitity in a.4}
\frac{\widetilde S}{H \ell} = \frac{\int_{0}^{\ell} (\ell - s)(s+d)^r \d s}{ \ell\int_{0}^{\ell}(s+d)^r  \d s},
\end{equation}
with $r = \frac{1}{p-1}$ and $d\geq0$. Now we obtain a lower bound for the ratio $\frac{\widetilde S}{H\ell}$. To this end, we define a probability measure $\mu_d$ by
\[
\mu_d(s) = \frac{(s+d)^r}{\int_{0}^{\ell} (\tau+d)^r \d \tau} ,\quad s\in[0,\ell],
\]
and also define its expected value by
$$
\mathbb{E}_{\mu_d }= \int_{0}^{\ell} s\mu_{d}(s) \d s.
$$
Then we may rewrite the ratio in \eqref{eqn: key idenitity in a.4} as
\[
\frac{\widetilde S}{H \ell} =1 -\frac{\mathbb{E}_{\mu_{d}}}{\ell}.
\]
Denote $\mu_* = \mu_0$. We now show that $\mathbb{E}_{\mu_*} > \mathbb{E}_{\mu_d}$ for all $d>0$. Since $\mu_*$ and $\mu_d$ are probability measures, we have
\begin{equation}
\label{eqn: int of mus1-mub=0}
 \int_{0}^{\ell}(\mu_*(s)-\mu_{d}(s))\d s = 0.
\end{equation}
Hence, since $0= \mu_*(0) < \mu_{d}(0)$, there must exist $s_c \in (0,\ell)$ such that $\mu_*(s_c) = \mu_d(s_c)$, as otherwise we would have $\int_{0}^{\ell} \mu_*(s) - \mu_{d}(s) \d s <0$, contradicting \eqref{eqn: int of mus1-mub=0}. By equating $\mu_*$ and $\mu_{d}$, we can find an explicit form for $s_c$ and deduce that it is unique. Let $M(s) = \int_{\ell}^s \mu_*(\tau)-\mu_{d}(\tau) \d \tau$. By \eqref{eqn: int of mus1-mub=0}, we have $M(\ell)=0 = M(0)$, and integration by parts yields
$$
    \mathbb{E}_{\mu_*} - \mathbb{E}_{\mu_d}=\int_0^\ell s M'(s)\d s =-\int _0^\ell M(s) \d s.
$$
Moreover,
$$
\operatorname{sgn}(M')(s) =\operatorname{sgn}(\mu_* - \mu_d)(s)= \begin{cases}
-1, & \text{for } s < s_c, \\
1,  & \text{for } s > s_c.
\end{cases}
$$
Hence, we deduce that $M(s)< 0 $ on $(0,\ell)$, which directly implies $\mathbb{E}_{\mu_*} > \mathbb{E}_{\mu_d}$ for all $d>0$. Then
\[
\frac{\widetilde S}{H\ell} = 1-\frac{\mathbb{E}_{\mu_d}}{\ell} \geq 1-\frac{\mathbb{E}_{\mu_*}}{\ell} = \frac{p-1}{2p-1}.
\]
Rearranging the terms yields
\begin{equation}\label{eqn:x-diff}
\ell \leq \frac{2p-1}{p-1} \frac{ \widetilde S }{H}
\end{equation}
Noting that $\widetilde S \leq S$ and combining this with \eqref{eqn:ineq in previous A.3} completes the proof when $f_0$ is of the form \eqref{eqn:Needyform2}.

\noindent\textbf{Step 3: The case where $f_0$ is of the form \eqref{eqn:Needyform1}.}
Since the argument is similar to the case where $f_0$ is of the form \eqref{eqn:Needyform2}, we only outline the main steps. By following the argument in Step 1 for \eqref{eqn:ineq in previous A.3}, we can prove
$$
\mathcal{F}_I^p(f_0) \geq \frac{H^p}{((s_4'-s_3')+(s_2'-s_1'))^{p-1}}.
$$
Then, since \eqref{eqn:x-diff} still holds, and following the argument in Step 2 used to prove  \eqref{eqn:x-diff}, we can additionally prove
$$
s_4'-s_3' \leq \frac{2p-1}{p-1} \frac{\widehat S}{H}, \quad \mbox{with }\widehat S=\int_{s'_3}^{s'_4} f_0(s)\d s .
$$
The result then directly follows by combining the two inequalities and noting $S_1+S_2 \leq S$. This concludes the proof.
\end{proof}

\bibliographystyle{abbrv}
\bibliography{References}

\end{document}